\documentclass{article}
\usepackage{amsfonts, amsmath, amsthm, amssymb}
\title{Conjugacy of Real Diffeomorphisms.\\
A Survey$^1$}
\author{Anthony G. O'Farrell$^2$\\
Mathematics Department\\
NUI, Maynooth, Co. Kildare\\
Ireland
\\
and
\\
Maria Roginskaya\\
Mathematics Department\\
Chalmers University of Technology and G\"oteborg University\\
SE-412 96 G\"oteborg\\
Sweden
}
\date{\today}

\def\bb#1{{\mathbb{#1}}}

\newtheorem{theorem}{Theorem}[section]%
\newtheorem{lemma}[theorem]{Lemma}%
\newtheorem{corollary}[theorem]{Corollary}%
\newtheorem{proposition}[theorem]{Proposition}%

\newtheorem{example}{Example}[section]

\def\Proof{{\par\medskip\noindent\bf Proof. }}
\def\qed{{\hfill\vrule height 4pt width 4pt depth 0pt
             \par\vskip\baselineskip}}
\def\Definition{{\par\medskip\noindent\bf Definition. }}

\def\implies{{\ \Rightarrow\ }}

\def\Oh{{\hbox{O}}}
\def\oh{{\hbox{o}}}
\def\clos{{\rm clos}}
\def\bdy{{\rm bdy}}

\def\half{{\raise1pt\hbox{$\scriptscriptstyle{1\over 2}$}}}
\def\ONE{1\hskip-4pt1}

\def\R{{\bb R}}
\def\Z{{\bb Z}}
\def\C{{\bb C}}
\def\N{{\bb N}}
\def\S{{\bb S}}

\def\Diffeo{{\textup{Diffeo}}}
\def\Homeo{{\textup{Homeo}}}
\def\Conj{{\textup{Conj}}}

\def\deg{{\textup{deg}}}
\def\fix{{\textup{fix}}}
\def\sign{{\textup{sign}}}

\begin{document}
\maketitle


\footnotetext[2]{Supported by Grant SFI RFP05/MAT0003 and the ESF Network
HCAA.}

\section*{Abstract}
Given a group $G$, the conjugacy problem in $G$ is the
problem of giving an effective procedure for determining
whether or not two given elements $f,g\in G$ are
conjugate, i.e. whether there exists $h\in G$
with $fh=hg$.  This paper is about the conjugacy problem
in the group $\Diffeo(I)$ of all
diffeomorphisms of an interval $I\subset\R$.

There is much classical work on the subject,
solving the conjugacy problem for special classes of maps.
Unfortunately, it is also true that many results and arguments
known to the experts are difficult to find in the literature,
or simply absent. We try to repair these lacunae,
by giving a systematic review,
and we also include new results about the
conjugacy classification in  the general case.

\section{Informal Introduction}
\subsection{Objective}
We are going to work with
diffeomorphisms defined
just on various intervals (open, closed, or half-open, bounded or
unbounded).
Let $\Diffeo(I)$ denote the group of (infinitely-differentiable)
diffeomorphisms of the interval $I\subset\R$, under  the
operation of composition.
We denote the (normal) subgroup of orientation-preserving
diffeomorphisms of the
interval by $\Diffeo^+(I)$.  If an endpoint $c$ belongs
to $I$, then statements about derivatives at $c$ should be interpreted
as referring to one-sided derivatives.

\medskip
Our objective is to classify the conjugacy classes,
i.e. to determine when
two given maps $f$ and $g$ are conjugate in $\Diffeo(I)$.

\medskip
The reason this problem is important is that conjugate elements
correspond to one another under a \lq\lq change of variables".
For most applications, a change of variables will not
alter anything essential, so only the conjugacy class
of an element is significant. From the viewpoint
of group theorists, it is also usual to regard only the conjugacy classes
as having \lq\lq real" meaning in a group.

\medskip
Throughout the paper, we will use the term {\em smooth}
to mean infinitely-differentiable. There is a good deal of
valuable and delicate work on conjugacy problems for functions that
are merely $C^k$, but we will not delve into this (apart from
an occasional remark), in order
to keep the discussion within bounds.

\medskip
Apart from its intrinsic interest,
the conjugacy problem has applications to
the holonomy theory of codimension-one foliations.
Mather established a connection between the homotopy of
Haefliger's classifying space for foliations and
the cohomology of
the group $G$ of compactly-supported diffeomorphisms of the line
\cite{M2, M3}.
Mather also used a conjugacy classification of a subset of
the group $G$
in order to establish that
$G$ is perfect. It follows from a result of Epstein
that $G$ is simple. Our own
study of the conjugacy problem arose independently from our
interest in {\em reversible maps}
(maps conjugate to their own inverses) \cite{OS}.

\medskip
We should make it clear that we are not here discussing an example
of the classical Dehn conjugacy problem of combinatorial
group theory. The group $\Diffeo^+(I)$ is not
countably-presented. It has the cardinality of the continuum.
Its family of conjugacy classes also has the cardinality of the
continuum.
To classify conjugacy classes is a matter of identifying suitable
conjugacy invariants which separate the classes. To be of practical
use, the invariants should be reasonably \lq\lq computable", in some
sense, but the sense has to be more lax than
standard Turing-machine computability.
For a start, we assume that we have available a \lq\lq real
computer", that can do real arithmetic and decide
equality of two suitably specified real numbers.
We include as suitable specifications things like
the value of an integral of a suitably-explicit function,
and the limit of a suitably-explicit sequence.
In practice, the kind of problem one wishes to solve is this:
given a prescription for two diffeomorphisms $f$ and $g$,
sufficiently explicit that we may compute the images
of any suitably-specified point, decide whether or
not they are conjugate. This may seem quite modest (especially
as we have not make explicit what is meant by
\lq\lq suitably-explicit")
but, as we shall see, it is rather too much to hope for.
A less demanding task would be to come up with a
procedure that will confirm that two non-conjugate
diffeomorphisms are in fact non-conjugate, but may
go on forever if presented with two conjugates.
Even this is too much, except in special cases.
What one can do is provide
a collection of classifying invariants
that provide a significant conceptual simplification of the conjugacy
problem.

\medskip
The conjugacy  problem
in $\Diffeo(I)$ may be reduced to the corresponding problem in the
subgroup $\Diffeo^+(I)$ of orientation-preserving maps of $I$ --- this recent
result is described in Section \ref{section-Diffeo}.
A crucial case of the latter problem
is the special case in which
the diffeomorphisms $f,g\in \Diffeo^+(I)$ are
fixed-point-free on the interior $J$ of $I$.
The problem is trivial if $I$ is open
(Proposition \ref{proposition-fpfree}).  A new result
(Theorem \ref{theorem-main})
provides an effective way to approach it when $I$ is
half-open. We establish that it suffices to
search for a conjugacy $h$ among the solutions
of a first-order ordinary differential equation.
This also helps with the case of compact $I$.
For special (\lq\lq flowable") diffeomorphisms of a compact $I$,
the conjugacy classification can be achieved using
the so-called \lq\lq functional moduli",
similar in character to the \'Ecalle-Voronin moduli for
the conjugacy classification of biholomorphic germs
\cite{V}. In the general case, this cannot be done.

There has been much work on this problem.
Important steps
in the story we describe below
are the work of Sternberg,
Takens, Sergeraert, Robbin, Mather,
Young, and
Kopell, among others. There is a useful
summary survey of progess up to
1995 by Ahern and Rosay \cite{AR}.  See also
references
\cite{S, T, SE, RO, M}, \cite[Chapter 8]{KCG},
\cite[Chapter 2]{KH}, \cite{MS, Y, B, ALY, N, SZ, YO}.
There are some parallels with the conjugacy problem
for complex analytic germs, for which see
\cite{CG}.

\subsection{Notation}
We shall use $\Diffeo$ as an abbreviation for $\Diffeo(I)$,
and $\Diffeo^+$ for $\Diffeo^+(I)$, whenever there is no
danger of confusion.

For $f\in\Diffeo(I)$, we denote the set of fixed points
of $f$ by fix$(f)$.

\medskip
We use the symbol $f^{\circ n}$ for the $n$-th
iterate of $f$ (i.e the $n$-th power in the
group $\Diffeo(I)$). We also use it for negative
$n=-m$, to denote the $m$-th iterate of
the inverse function $f^{\circ-1}$.
The notation
$f^{\circ0}$
denotes
the identity  map $\ONE$.

\medskip
We use similar notation for compositional powers and inverses
in the group $F$
of formally-invertible formal power series (with real coefficients)
in the indeterminate
$X$. The identity $X+0X^2+0X^3+\cdots$ is denoted simply by $X$.

\medskip
\noindent We denote
$g^h=h^{\circ-1}\circ g\circ h$, whenever
$g,h\in\Diffeo(I)$.   We say that {\em $h$ conjugates $f$ to $g$}
if $f=g^h$.

\medskip
We use the notation deg$f$ for the degree of the diffeomorphism
$f\in\Diffeo(I)$ ($=\pm1$, depending on whether or not
$f$ preserves the order on $I$).
Thus
$$ \Diffeo^+(I) =\{f\in\Diffeo(I): \deg f = +1\}.$$

Given a closed set $E\subset I$,
we set
$$ \Diffeo^+_E = \Diffeo^+_E(I) = \{f\in\Diffeo^+(I): f(x)=x,\ \forall x\in E\},$$
the subgroup of those direction-preserving maps that fix
each point of $E$.

We denote the map $x\mapsto -x$ on $\R$ by $-$.

When we come to discuss conjugation for elements
having complicated fixed-point sets, we will need
notation for the available conjugacies on
particular subintervals. So we make a definition:

\Definition Let $f,g\in \Diffeo^+(I)$. Given an open interval
$J\subset I$
that is mapped onto itself
by $f$ and by $g$,
we say that a map $\phi\in\Diffeo(\clos(J))$
is a $J$-{\em conjugation} from
$f$ to $g$ if $f^\phi=g$ on $\clos(J)$.  We denote the set of all
$J$-conjugations from $f$ to $g$ by
$\Conj(f,g;J)$, or just $\Conj(J)$, if the context is clear.

\subsection{Remarks about Topological Conjugacy}
A necessary condition for the conjugacy of two elements
$f,g\in\Diffeo(I)$ is that they be topologically-conjugate,
i.e. conjugate in the homeomorphism group $\Homeo(I)$.

The homeomorphism problem is strictly easier than the diffeomorphism
problem, because it is included as part of it: (1) One can
show that each conjugacy class of homeomorphisms has an
element that is a diffeomorphism. (2) Thus, if one knows how
to classify diffeomorphisms up to topological conjugacy,
then one knows how to classify homeomorphisms also. (3)
the topological conjugacy classification is
coarser
than the diffeomorphic.

As we shall now explain, the topological conjugacy problem
is already intractable, in computational terms, so it follows
that the same is true for smooth conjugacy.

Let us consider the case $I=\R$.

The conjugacy problem in the homeomorphism group $\Homeo(\R)$
has a classical solution in terms of a \lq\lq symbol" invariant.
This goes back, essentially, to Sternberg \cite{S}, who in 1957 described
the conjugacy classes in the group of germs of homeomorphisms of
neighbourhoods of a point on the line. For an exposition of the
classification in $\Homeo(\R)$, see \cite{OF1}.
Mere
topological conjugacy of two direction-preserving diffeomorphisms
$f$ and $g$ is determined by the existence of a homeomorphism of
$\R$ mapping $\fix(f)$ onto $\fix(g)$ and coincidence of the \lq\lq
pattern of signs" of $f(x)-x$ and $g(x)-x$ off the fixed-point sets.
(The pattern of signs of $f(x)-x$ is called the \lq\lq signature" of
$f$.)

Suppose $f=g^h$. Then $h$ carries $F_1=\fix(f)$ onto $F_2=\fix(g)$, so the pairs
$(\R,\fix(f))$ and $(\R,\fix(g))$ are homeomorphic.
An {\em order isomorphism} between two partially-ordered sets
is an order-preserving bijection. Two partially-ordered sets
are {\em order-isomorphic} if there exists an order isomorphism between
them. Order-isomorphism is an equivalence relation on the
family of partially-ordered sets, and the equivalence classes
are called {\em order classes}.
The homeomorphism class of a pair
$(\R,F)$ (with $F$ closed) is determined by the order class
of $F$ (with the usual
total
order inherited from $\R$).
Every closed subset of $\R$ is the fixed-point set of some
homeomorphism (and even of some diffeomorphism), so
the set $\fix(f)$ may be quite general.
Thus there are two obstacles to finding an algorithmic solution
to the topological conjugacy problem:

\medskip
(1)
The problem
of determining whether two closed subsets of $\R$ are order-equivalent
does not appear to be amenable to an algorithmic solution.
For subsets of simple structure it may be resolved by
noting that an order-isomorphism will induce a bijection of
the derived set, the second derived set, and so on through ordinals,
a bijection of the relative complements of each of these,
a bijection of the condensation set,
a bijection of each interval subset, and of the ends of
such intervals, and of derived sets of ends, etc.
But a general algorithm is another matter.

\medskip
(2) There
may be a large collection of order-isomorphisms
between $\fix(f)$ and $\fix(g)$, and we then
need some systematic way to check for the existence
of one order-isomorphism that gives
a concidence of signatures on the corresponding
complementary intervals.

Returning to the problem of $C^\infty$ conjugacy,
we have additional complications,
as the following observations indicate.

\subsection{Smooth Conjugacy of Pairs}
Suppose $f,g,h\in\Diffeo(I)$ and
$f=g^h$. Then $h$ carries $F_1=\fix(f)$ onto $F_2=\fix(g)$, so the pairs
$(\R,\fix(f))$ and $(\R,\fix(g))$ are diffeomorphic. This necessary
condition is more complex to check than the corresponding
topological condition.
To determine
whether two homeomorphic pairs belong to the same diffeomorphism class, it
is necessary to search among all the order-isomorphisms
of the $F_i$'s for one having a
diffeomorphic extension. The existence of a diffeomorphic
extension may be checked using a theorem of
Whitney.
Whitney's condition \cite{W} for the existence of a $C^\infty$
extension of a function $h$ from $F_1$ to $\R$ may be stated as
follows: For each $k\in\N$, the $k$-th Newton divided
difference of $h$ is uniformly continuous on bounded
subsets of
$$ \{(x_1,\ldots,x_{k+1})\in
F_1\times\cdots\times F_1:
x_i\not=x_j,\ \forall i\not= j\}$$
(i.e. extends continuously to the full product
$F_1\times\cdots\times F_1$).
In fact, an order-isomorphism has a diffeomorphic extension
if and only if it has an infinitely-differentiable extension,
and the (uniquely-determined) first derivative of such
an extension is nonzero at each accumulation point.

\subsection{Orbits, Multipliers and Taylor Series}
The (two-sided) {\em orbit} of a point $a\in\R$ under a diffeomorphism $f$
is the set $\{f^{\circ n}(a): n\in\Z\}$ of all forward
and backward images of  $a$ under the action of $f$.
If $f,g,h\in\Diffeo(\R)$ and $f=g^h$, then for each $a\in\R $, the map $h$ carries
the orbit $O_1$ of $a$ under $f$ onto a corresponding orbit
$O_2$ under $g$, so the pairs $(\R,O_1)$ and
$(\R,O_2)$ are diffeomorphic.  An implication is that these pairs
are equivalent under locally-bi-Lipschitzian maps. Thus, for instance,
one sees (by estimating the number of points in
orbits in intervals of comparable length) that the maps
defined by $$ f(x) = x + \exp(-1/x^2) $$
and
$$ g(x) = x + \exp(-2/x^2) $$
are not conjugate, although they have identical signatures.

It is straighforward (using Whitney's result) to check whether
two orbit-pairs $(\R,O_1)$ and $(\R,O_2)$ are
diffeomorphic, but a difficulty is that one must check that
for each orbit of $f$ {\em there exists} some orbit
of $g$ that gives a diffeomorphic pair. This is not a constructive
condition, as it stands.

\medskip
Further, if $f=g^h$, then for each $a\in\fix(f)$, letting
$b=h(a)$, we have
$$ g^\prime(b) = f^\prime(a), $$
i.e. $f$ and $g$ have the same \lq\lq multipliers" at corresponding points.
This necessary condition actually follows from
the previous one about orbits when $a$ is a boundary
point of $\fix(f)$, but is easier to check
when it fails. It is trivial at accumulation points
of $\fix(f)$.

There is a more elaborate
necessary condition involving higher
derivatives, best expressed in terms of Taylor
series:
Let $T_af$ denote the truncated Taylor series of $f$ about $a$:
$$ T_af = \sum_{n=1}^\infty \frac{f^{(n)}(a)}{n!} X^n $$
(regarded as a formal power series in an indeterminate $X$).
One then has
$$ T_af = (T_ah)^{\circ-1}\circ (T_bg)\circ (T_ah),$$
where $\circ$ denotes the formal composition, and
$p^{\circ-1}$ denotes the formal compositional inverse.
This condition is weaker than the one about orbits,
since the Taylor series at a point is determined
by the values of the function at any sequence tending to
the point.

However, there is a straightforward algorithm for checking whether
or not two formal power series are formally conjugate.
In fact, each series is conjugate to one of $\lambda X$ ($\lambda\in \R$),
or one of $\pm X\pm X^{p+1}+\alpha X^{2p+1}$, and in
each case the correct class can be determined by a terminating computation.
This fact is well-known (cf. \cite[p.546]{AR}, \cite{L, Ka, OF2}), and is routine
to check.  The main point to note is that the group
of invertible formal power series (with its product topology)
is topologically-generated by
the maps $x\mapsto\lambda X$ ($\lambda\not=0$) and
$x\mapsto x +\alpha x^{p+1}$ ($p\in\N$)\cite[Lemma 1, p.5]{OF2}.

For example:
\begin{enumerate}
\item $3X+ X^2$ is conjugate to $3X+2X^2$, and to any other series
that begins with $3X$, but is not conjugate to any series that
begins with $2X$;
\item $X+ X^2 + X^3$ is conjugate to
$X+2X^2+4X^3+8X^4+\cdots$, but not to any series beginning
with $X+3X^2+6X^3$ or $X+2X^3$;
\item Each series beginning $X+X^4+2X^7$  is conjugate to
each series beginning $X+5X^4+50X^7$.
\end{enumerate}

We will see below that there is more to conjugacy than
the diffeomorphism of pairs, correspondence of signatures,
and the orbit conditions,
but that the problem can nevertheless be reduced
to manageable proportions, provided one does not try to
do the impossible.

\subsection{Centralisers}
Typically, if $f$ and $g$ are conjugate
diffeomorphisms, then the family $\Phi$
of diffeomorphisms $\phi$ such that
$f=\phi^{\circ-1}\circ g\circ\phi$
has more than one element. In fact
$\Phi$ is a left coset of the centraliser
$C_f$ of $f$ (and a right coset of $C_g$).
For this reason, it is important for us to
understand the structure of these centralisers.
The problem of describing $C_f$ is a special
conjugacy problem --- which maps conjugate
$f$ to itself?

Historically, there has been a good deal more
work on the problem of centralisers than on the
general conjugacy problem.

There may be a great many conjugacies between two given
conjugate diffeomorphisms. In the open-interval case, the
centraliser of
a fixed-point-free diffeomorphism is very large, and is not
abelian.

Kopell \cite{K} showed
that when $I$ has one of its endpoints as a member,
then the centraliser
of an $f$ that is fixed-point-free on the interior
of $I$ must be quite small --- it is a subgroup
of a one-parameter abelian group, and it may consist
just of the iterates of $f$.
An example was
given by Sergeraert \cite{SE};  probably this behaviour
is \lq\lq generic".
Kopell \cite{K} showed that it is
generic for maps that fix
only the two endpoints of $I$.
These phenomona tell us that in many cases the search
for a conjugating map $h$ from $f$ to $g$ may be
confined to a 1-parameter search space. Our main new theorem
gives a specific way to locate this search space,
in the case of a half-open interval.

See Subsection \ref{subsection-flowability}
and Section \ref{section-sufficient}).

\subsection{Outline}
The paper is organised as follows.

First, we consider various special cases
of the full conjugacy problem, and related
simpler problems, and then we use these cases and
problems as building blocks in constructing a
solution to the full problem.

The results are summarised formally in Section \ref{section-statements}.
The remaining sections provide proofs, elaboration, and examples.

In less formal terms, we proceed as follows:

We start with
the simple and classical case of
fixed-point-free maps of an open interval,
where there is just one conjugacy class
of diffeomorphisms.
(Details are in Section \ref{section-fpfree}.)

Then we study conjugacy in half-open intervals,
starting with diffeomorphisms of the interval $[0,+\infty)$
that fix only $0$.
First, we review the classical results based on
normal forms that exist when the diffeomorphism is not tangent to
the identity to infinite order.
In these special cases the conditions simplify. If the multiplier at
$0$ is not $1$ (i.e. $0$ is a {\em hyperbolic} fixed point), then
Sternberg \cite{S} identified the multiplier as the sole conjugacy
invariant. If the multiplier is $1$, but $f$ is not \lq\lq
infinitesimally tangent to the identity" (i.e. $T_0f\not=X$ --- we
find it less of a mouthful to express this condition as \lq\lq $f-x$
is not flat at $p$"), then Takens \cite{T}
identified the conjugacy class
of the Taylor series $T_0f$ (in the group
of formally-invertible power series) as the sole conjugacy invariant.
We show that the general problem cannot be tackled
using normal forms. We identify an infinite product condition that is
necessary for conjugacy. We then
base our approach to characterising
conjugacy on a certain differential equation that may be formulated
when the product condition holds.
(Details are in Sections \ref{section-half-open}, \ref{section-necessary}
and
\ref{section-sufficient}.)

Next, we study
conjugacy in $\Diffeo^+(I)$, for
closed bounded intervals $I$, for maps that are
fixed-point-free on the interior $J$ of $I$.
In the \lq\lq Axiom A" case, in which
both fixed points are hyperbolic, Robbin
characterised conjugacy in terms of the multipliers
and a \lq\lq modulus" (a smooth function on
$(0,+\infty)$; detail below).  Results
of Young \cite{Y} relate to other cases in which
$f-x$ is not flat at either end of $I$,
particularly the \lq\lq saddle-node" case,
in which $T_0f-X$ is zero mod$X^2$, but
not zero mod$X^3$, for both endpoints $p$.
He used  so-called \lq\lq formal multipliers"
(certain diffeomorphisms from $J$ to $(0,+\infty)$)
to construct a substitute for the Robbin modulus,
which, when taken
together with the conjugacy classes of the
Taylor series at the ends, characterise conjugacy
classes. There is a more general treatment of
functional moduli ideas in unpublished work of Mather \cite{M2}.
We give a necessary and sufficient condition
for conjugacy that builds on the result for half-open intervals.
We also review functional moduli in the special Mather
case, and a useful new necessary condition
expressed in terms of the \lq\lq shape"
of a graph associated to the pair of maps $(f,g)$.
(Details are in Section \ref{section-compact}.)

Then we move on to general direction-preserving diffeomorphisms,
on any interval $I$,
with possibly complicated
fixed-point sets. We take this in two stages:
\\
(1)
We reduce
the conjugacy problem in $\Diffeo(I)^+$to
the conjugacy problem in $\Diffeo^+_{\bdy E}(I)$, for a fixed closed $E$.
(Details in Subection \ref{subsection-Diffeo+}).
\\
(2)
We address the conjugacy
problem in $\Diffeo^+_{\bdy E}(I)$ for maps that
belong to $\Diffeo^+_E(I)$ and are fixed-point-free off $E$.
(Details in Subsection \ref{subsection-Diffeo+B})

The final theoretical step is
the reduction of the conjugacy problem in $\Diffeo(I)$
to the conjugacy problem in $\Diffeo(I)^+$.
(Details in Section \ref{section-Diffeo}).

\medskip
By a {\em flow} on an interval $I$, we mean
a continuous homomorphism $t\mapsto \Phi^t$ from the additive
topological group
$(\R,+)$ into $\Diffeo^+(I)$, endowed with its usual topology
(the topology of simultaneous convergence of functions and their
inverses,
together with all their derivatives,
uniformly on $I$).

We say that $f\in \Diffeo^+(I)$ is {\em flowable} if there
exists a flow $\Phi^t$, with $f=\Phi^1$ (i.e.
$f$ is the \lq\lq time 1"
map of the flow $(\Phi^t)_{t\in\R}$.

There is a connection between
our subject and the question of when an $f\in\Diffeo^+(I)$ is flowable.
For this, see also \cite{SE}. We shall make some remarks
about flowability as we go along (cf. Subsection
\ref{subsection-flowability})

\medskip
Along the way, we present some conjectures and problems that,
if true or solved,
as the case may be,
would improve our understanding of one-dimensional conjugacy.

\section{Overview and Statement of Main Results}
\label{section-statements}
\subsection{Open Intervals}
A fixed-point-free diffeomorphism of an open interval $I$
 must preserve orientation.
There is just one conjugacy class of fixed-point-free
diffeomorphism in $\Diffeo(I)$,
which splits into just two conjugacy classes
with respect to $\Diffeo^+(I)$:

\begin{proposition}[Sternberg \cite{S}]\label{proposition-fpfree}
Suppose $I$ is an open interval and
$f$ and $g$ are fixed-point-free
elements  of $\Diffeo(I)$.
Then $f$ and $g$ are conjugate in
$\Diffeo^+(I)$ if and only if their
graphs lie on the same side of the diagonal.
\qed
\end{proposition}

This is proved in Section \ref{section-fpfree} below.

\subsection{The Interval $[0,+\infty)$}
Note that $\Diffeo(I)=\Diffeo^+(I)$ whenever $I$
is a half-open interval, because all the elements
of $\Diffeo(I)$ have to fix the endpoint that belongs
to the interval.

Consider $f,g\in \Diffeo([0,\infty))$, fixed-point-free on
$(0,\infty)$. Under what
circumstances does there exist an $h\in \Diffeo([0,\infty))$ with
$f=g^h$?

The set of all $f\in \Diffeo([0,\infty))$, that fix only $0$ is
the disjoint union of the two subsets
$$S_+=\{f:f(x)>x \mbox{  on  }(0,\infty)\}$$
$$S_-=\{f:f(x)<x \mbox{  on  }(0,\infty)\}$$
each of which is a sub-semigroup of $\Diffeo([0,\infty))$. Each of
these semigroups is preserved by conjugacy, i.e. is a union of
conjugacy classes. Thus, for $f$ to be conjugated to $g$ it is
necessary that {\sl they belong to the same semigroup, $S_+$ or $S_-$}.
We call this the {\bf \lq\lq the sign condition".}

\medskip
{\bf Remark.}
In later sections, where the context changes,
the meaning of \lq\lq the sign condition" will
change as well.  So the above defines the sign condition
just for the case of $\Diffeo([0,+\infty))$.

\medskip
Note that $f\in S_+$ is equivalent to $f^{\circ -1}\in S_-$, so that
to characterize conjugacy it suffices to consider $f\in S_-$.

We review some special cases, and then look at the
general case.
\subsection{$[0,+\infty)$: Hyperbolic Case}
The result for the case $f'(0)\not=1$ is
known as Sternberg's Linearization Theorem.
It was essentially proved in \cite{S}.
It may be regarded as the smooth equivalent
of Schroeder's theorem \cite[Chapter II]{CG} about complex analytic germs
in one variable.
\begin{theorem}\label{notone} Let $f,g\in S_-$ and $f'(0)\neq 1$. Then the following are equivalent:
\begin{enumerate}
    \item $f'(0)=g'(0)$;
    \item There exists $h\in \Diffeo^+([0,\infty))$ with $f=g^h$;
    \item For each $\lambda>0$ the sequence $h_n=g^{\circ - n}(\lambda f^{\circ n})$
    converges (pointwise)
    to a diffeomorphism $h$ on $[0,\infty)$;
    \item The sequence $h_n=g^{\circ - n}\circ f^{\circ n}$ converges to a diffeomorphism
    $h$ on $[0,\infty)$;
    \item There exists $\lambda>0$, such that the sequence $h_n=g^{\circ - n}(\lambda f^{\circ n})$
    converges to a diffeomorphism $h$ on $[0,\infty)$;
\end{enumerate}
\end{theorem}
The details are in Subsection \ref{SLT}.

\begin{corollary}\label{corollary-flow-Sternberg}
If $f\in S$, and $f'(0)\not=1$, then the centraliser
$C_f$ of $f$ in $\Diffeo([0,+\infty))$
is a one-parameter group, and $f$ is flowable.
\end{corollary}

\Proof  In fact, $C_f$ is the image under an inner
automorphism of
\\ $\Diffeo([0,+\infty))^+$ of the centraliser
of $x\mapsto f'(0)x$, and this consists precisely
of the maps $x\mapsto\mu x$ with $\mu>0$.

Also, the map $x\mapsto f'(0)$ is the time 1
map of the flow $t\mapsto f_t$, where
$$f_t(x) = e^{t\ln f'(0)}x.$$
\qed

\subsection{$[0,+\infty)$: Taylor Series}

Since all the elements fix $0$,
we see $f=g^h$ in $\Diffeo^+([0,+\infty))$  implies
$$ T_0f = (T_0h)^{\circ-1}\circ (T_0g)\circ (T_0h),$$
Thus $T_0f$ and $T_0g$ are conjugate in the group
of formally-invertible series.  We call this
{\bf Condition (T).}

In case $f'(0)\not=1$, condition (T)
just says $f'(0)=g'(0)$.  In the
non-hyperbolic case,
it imposes conditions on some higher derivatives.

For the non-flat case, Takens
\cite[Theorem 2]{T} proved the following theorem\footnote{
There is folklore that says
that Mather
independently found this result, but
we checked with Mather, who said he definitely did not.}.

\begin{theorem}[Takens]\label{theorem-2.4}
Suppose that $f,g\in S_-$,
and $f-x$ is not flat at $0$.
Then the following are equivalent.
\\
(1) Condition (T) holds.
\\
(2) There exists $h\in\Diffeo([0,+\infty))$ such that
$f=g^h$.
\end{theorem}

Note that this generalises the equivalence of (1) and (2)
in Theorem \ref{notone}, since the multiplier
determines the conjugacy class of the series when it is not $1$.

See Subsection \ref{subsection-takens} for detail.

\subsection{$[0,+\infty)$: The Case $f-x$ flat at $0$}
If $f-x$ is flat at $0$, Condition (T)
just says that $g-x$ is also flat at $0$.
This is not enough.

\begin{example} Let $f(x)=x-e^{-\frac1{x}}$ and $g(x)=x-e^{-\frac1{x^2}}$.
The functions $f$ and $g$ are not conjugate in $\Diffeo([0,+\infty))$.
\end{example}
\Proof Suppose $h\in\Diffeo([0,+\infty))$, with
Taylor series $T_0h=aX+bX^2+\ldots$, is a conjugation.
Then it maps the interval $\left[\displaystyle\frac{x}2,x\right]$
to the interval
$\left[\displaystyle\frac{ax}2+o(x),ax+o(x)\right]$.
For small positive $x$, the first
interval has no more than $x\exp(2/x)$ iterates of $x$ under $f$,
whereas the second has at least
$(x/2)\exp(1/4a^2x^2))$
iterations
of $h(x)$ under $g$, a much greater
number.
But
the conjugacy condition
requires that the two intervals contain equal numbers of
iterates of $x$ and $h(x)$, respectively.
\qed

So we need another idea, in order to deal with
two general elements $f,g\in S_-$.
If you think about it, the main difficulty of the conjugacy problem
of the present section involves the functions
with $f-x$ flat at $0$. When endowed with the
relative topology from the usual Frechet-space
topology on $C^\infty([0,+\infty)$,
the group $\Diffeo^+([0,+\infty))$
is separable and metrisable, so has the cardinality
of the continuum, and hence (since Sternberg
gives us a continuum of conjugacy classes) the family
of conjugacy classes has the same cardinality.
From this point-of-view, Sternberg's family
is a substantial family of conjugacy classes.

However, cardinality is very crude way to measure
size. Another way is to use dimension.  The map
$D:f\mapsto f'$ is a continuous bijection between
$\Diffeo([0,+\infty))$ and the cone
of all smooth positive functions $h$ on $[0,+\infty)$
that have $$\int_0^\infty h(x)\,dx=+\infty.$$
This gives a way to embed our group as a convex subset
in a
Frechet space, and talk about the linear dimension
and codimension of subvarieties.  Sternberg's family
is the complement of a codimension $1$
variety, and so is a large part of the group.
But consider the conjugacy classes.
Conjugacy does not respect the convex
structure of the cone (i.e. it does not commute
with convex combinations), so we cannot induce
a manifold structure on the conjugacy classes.
What we can do is measure the size of a family
$\mathcal F$ of conjugacy classes in terms of the
minimal dimension of $R$, where $R$ ranges over
varieties that have at least one representative
of each  element of $\mathcal F$.  Let's call this
cardinal the {\em conjugacy dimension} of the family.
From this point of view,
Sternberg's family has conjugacy dimension 1.

If we take
$G_0=\Diffeo^+([0,+\infty))$
and denote by $G_n$ the subgroup consisting
of those $f\in G_0$
such that $f-x=$o$(x^n)$ at $0$, then we have
a countable nested chain of closed normal subgroups
$$ \cdots G_n\subset G_{n-1}\subset\cdots
\subset G_2\subset G_1\subset G_0 $$ and each difference set
$G_n\sim G_{n-1}$ has a naturally-parametrised one-parameter family
of conjugacy classes, so has conjugacy dimension one. Moreover, each
difference is an open dense subset of the next group $G_{n-1}$, so
it looks as though we have a nice stratification of the conjugacy
classes, with just a trivial collection left at the core. But the
fun really starts when we move inside the intersection $G_\infty$ of
the chain. For instance, to each element $\phi\in G_\infty\cap S_+$
we may associate a normal subgroup
$$ G_\phi = \{f\in G_0:
f(x)-x = O(\phi(x) \}.$$
One sees that the intersection of each countable
family of groups $G_\phi$ is nontrivial, so by
transfinite induction one can construct
uncountable chains of $G_\phi$'s.
From the purely algebraic point-of-view,
this is no different from what one can do
inside the Sternberg family, because one can construct
uncountable chains of normal subgroups
by restricting the multiplier to
subfields of the reals.  But from the analytic
point-of-view the $G_\phi$ are quite
different groups, because their images
under $D$ are cones, and invariant under
multiplication by positive reals.
This makes it clear
that there is no hope of tackling the conjugacy
problem by reducing to explicit normal forms,
since the set-theoretic difference of two
normal subgroups is a union of conjugacy classes.

Neither is it possible to reduce it
to the temptingly straightforward task of
comparing vectorfields whose exponentials are the
given functions, for the simple reason that the
exponential map is not surjective [SE].
It is easy to check if two
flows are related by a smooth change of variables,
but not all diffeomorphisms are flowable.

{\em The only way to come at it is to take two functions
and compare them directly with one another, rather
than with some collection of templates.}

We find such a procedure by using a suitable infinite product, and differential
equation.

Arising from this discussion, we state a conjecture:

\medskip\noindent
{\bf Conjecture.} {\sl The conjugacy dimension of
the diffeomorphism group of $[0,+\infty)$
is uncountable.}

\subsection{$[0,+\infty)$: The Product}
Let us begin again, with general $f,g\in S_-$.
For $x>0$ and $\xi>0$, let
\begin{equation}\label{equation-basic-product}
H_1(x,\xi) = H_1(f,g;x,\xi)
=\prod\limits_{n=0}^{\infty}\frac{f'(f^{\circ n}(x))}{g'(g^{\circ n}(\xi))}
.\end{equation}
We say that $f$ and $g$ satisfy {\bf Condition (P)} if there exist
$x>0$ and $\xi>0$ such that the product $H_1(x,\xi)$
converges (to a nonzero limit).

The product $H_1(x,\xi)$ appears already in Sternberg's paper
\cite{S}, in the special case $g(x)=\lambda x$, and in Kopell's
paper \cite{K} in the case $f=g$. We have not seen it used
in the literature for general $f$ and $g$.

We shall show (Corollary \ref{corollary-one-all})
that if Condition (P) holds, then $H_1(x,\xi)$
exists for all $x>0$ and $\xi>0$, and
(Lemma \ref{lemma-smooth-products})
is infinitely-differentiable
and positive.  We may then consider the three-parameter
initial-value problem
\begin{equation} \label{equation-IVP-D1}
D_1(a,\alpha,\lambda):\qquad
\left\{
\begin{array}{rcl}
 \displaystyle \frac{d\phi}{dx} &=&
H_1(x, \phi(x))\lambda,\ \forall x>0,\\
\phi(a)&=&\alpha
\end{array}
\right.
\end{equation}
depending on $\lambda>0$, $a>0$ and $\alpha>0$.
We shall show
that for each given $a>0$ and
$\alpha>0$, there exists (Lemma \ref{lemma-exactly-one})
exactly one $\lambda>0$ for
which the (unique) solution $\phi$
to problem $D_1(a,\alpha,\lambda)$ has $f(a)=g^\phi(a)$,
and
(Lemma \ref{lemma-conjugation}) that this $\phi$ conjugates
$f$ to $g$ in $\Diffeo((0,+\infty))$,
and (Lemma \ref{lemma-one-sided})
extends in $C^1([0,+\infty))$, with $\phi'(0)=\lambda$.
We denote this unique $\lambda$ by $\Lambda_+(a,\alpha)$,
and the unique $\phi$ by $\Phi_+(a,\alpha)$.

Thus, subject to Condition (P), there is a
1-parameter family of $C^1$ conjugations from $f$ to $g$ on
$[0,+\infty)$\footnote{$\{\Phi_+(a,\alpha):a>0,\alpha>0\}$
is a 1-parameter family, because
$$\Phi_+(a,\alpha) = \Phi_+(b,\Phi_+(a,\alpha)(b))$$
for each $b>0$.}.
This immediately gives us a result about $C^\infty$
conjugacy on $[0,+\infty)$:
\begin{theorem}[Main Theorem]\label{theorem-main}
 Let $f,g\in S_-$.  Then $f$ is conjugate to $g$
in $\Diffeo([0,+\infty)$ if and only if Condition (P) holds and
there exists some $a>0$ and $\alpha>0$ for which $\Phi_+(a,\alpha)$
is $C^\infty$ at $0$.
\end{theorem}
The value of this result is that it narrows the search for
a conjugating map $\phi$ to the 1-parameter family of solutions
of an explicit ordinary differential equation.

\medskip
We repeat (for emphasis) the fact already noted
that when $f-x$ is not flat at
$0$, then Condition (T) implies $f$ is conjugate to $g$.
Thus, since Condition (T) is easier to check
than Condition (P),
the theorem is only interesting when $f-x$ is flat at $0$.

\subsection{General Half-open Intervals}
All the above results about $[0,+\infty)$
carry over to
diffeomorphisms of an arbitrary half-open
interval that
fix only the endpoint that belongs to the interval.
Each such interval is diffeomorphic to
$[0,+\infty)$.

For a general half-open interval $I=[d,c)$ or $I=(c,d]$,
we take $J=$int$(I)$ and
define $S_-$ as the semigroup of diffeomorphisms
$f\in\Diffeo(I)$ which iterate all points of $J$ towards the endpoint
$d$, and $S_+$ as the semigroup of those that
iterate all points of $J$ towards $c$.
In order to adapt the above results about
$f,g\in S_-$ to the interval
$J\cup \{d\}$, one should replace $(0,+\infty)$ by $J$,
and $0$ by $d$.
Then, for $f,g\in S_-$, the product condition (P)
takes precisely the same form (\ref{equation-basic-product}),
and the differential equation also, except that its domain is
the interior $J$.  The theorem yields,
by conjugating $I$ to $[0,+\infty)$, a precisely similar
result for $f,g\in S_-$ on $I$.

For future reference, we formulate the
condition (for two $f,g\in S_-$ satisfying condition (P)):

\noindent{\bf Condition (E):}\\
{\sl There exists
$a, \alpha\in J$,
for which the $C^1$ extension of
the solution
$\phi=\Phi_+(a,\alpha;\bullet)$
from $J$ to the point $d$
is actually $C^\infty$.}

\medskip
\noindent
It is equivalent to
replace  \lq\lq there exist $a,\alpha$" by
\lq\lq for each $a$ there exists $\alpha$".

In these terms,  we may state:

\begin{corollary}\label{corollary-half-open-}
Let $I$ be half-open. Two elements $f,g\in S_-$ are conjugate
in $\Diffeo(I)$ if and only if they satisfy
conditions (P) and (E).
\end{corollary}
\begin{proof}
Immediate from Theorem \ref{theorem-main}.
\end{proof}

For $f,g\in S_+$,
one applies this result to $f^{\circ-1}$ and
$g^{\circ-1}$, which lie in $S_-$.  Unwinding the definitions,
we see that {\bf Condition (P) for elements of $S_+$}
involves the infinite product
\begin{equation}
H_2(x,\xi)=
\prod_{n=1}^\infty{
\frac{g^\prime(g^{\circ-n}(\xi))}%
{f^\prime(f^{\circ-n}(x))}
},
\end{equation}
(for $x,\xi\in J$)
and the differential equation takes the
form:
\begin{equation}
\label{equation-IVP-D2}
D_2(a,\alpha,\mu):\qquad
\left\{
\begin{array}{rcl}
\displaystyle \frac{d\phi}{dx} &=&
{H_2(x, \phi(x))\mu}, \hbox{ on }J\\
\phi(a)&=&\alpha,
\end{array}
\right.
\end{equation}
for $a,\alpha\in J$ and $\mu>0$.

Assuming condition (P), one has, for each
for $a,\alpha\in J$, the existence of a
unique $\mu>0$ (denoted $\Lambda_-(a,\alpha)$)
for which the unique solution $\phi=\Phi_-(a,\alpha)$
has a $C^1$ extension to $d$, with $\phi'(d)=\mu$.
The version of {\bf Condition (E) for elements
of $S_+$} then says:

\medskip
{\sl There exists
$a, \alpha\in J$,
for which the $C^1$ extension of
the solution
$\phi=\Phi_-(a,\alpha;\bullet)$
from $J$ to the point $d$
is actually $C^\infty$.}

\medskip
With this terminology,  the previous corollary yields:

\begin{corollary}\label{corollary-half-open+}
Let $I$ be half-open. Two elements $f,g\in S_+$ are conjugate
in $\Diffeo(I)$ if and only if they satisfy
the
$S_+$ versions of  conditions (P) and (E).
\end{corollary}
\qed

\subsection{Compact Intervals}\label{subsection-compact}

Now we consider the question of conjugacy for
orientation-preserving diffeomorphisms
of a {\em compact} interval $I$, which
are fixed-point-free on the interior
$J$.

Let $f$ and $g$ be two such diffeomorphisms.

The first necessary condition is the same as before:

\medskip
\noindent
{\bf The sign condition:}
$ \sign(f(x)-x) = \sign(g(x)-x),\ \forall x\in J$.

\medskip
This means that $f$ and $g$ are topologically conjugate,
and have similar dynamics on $I$.  Forward iteration from any point
of $J$ converges monotonically to one end of $I$, and backward iteration leads to the
other end. So
the map $f$ induces a direction on $J$ ---
upward if $f(x)>x$ on $J$, downward if $f(x)<x$
on $J$. We label the ends of $J$
accordingly:
$$ d = d(J) = \lim_{n\rightarrow \infty} f^{\circ n},$$
$$ c = c(J) = \lim_{n\rightarrow\infty} f^{\circ -n}.$$
We call $c$ the \lq\lq initial endpoint" of $J$, and
$d$ its \lq\lq final endpoint". We call the direction towards
$d$ the \lq\lq forward direction" on $J$, and the
other the \lq\lq backward direction".

For a compact interval $I=[c,d]$, with nonempty interior $J$,
we define $S_-(I)$ as the semigroup of homeomorphisms
that iterate each element of $J$ towards $d$.

In order that two given $f,g\in S_-$ be conjugate in
$\Diffeo([c,d])$, it is necessary that they be
conjugate in
$\Diffeo([c,d))$ and in $\Diffeo((c,d])$). Thus Corollaries
\ref{corollary-half-open-}
and \ref{corollary-half-open+}
apply, and  tell us that
the two-sided product
\begin{equation}\label{expression-two-sided-product}
H(x,\xi) = H(f,g,x,\xi)=H_1(x,\xi)/ H_2(x,\xi) =
\prod_{n=-\infty}^\infty
\frac{f^\prime(f^{\circ n}(x))}%
{g^\prime(g^{\circ n}(\xi))}
\end{equation}
must converge for some (or equivalently all) $x,\xi\in J$.
This is the appropriate version of {\bf Condition (P),
for compact intervals}.

Assuming Condition (P), we may form two initial-value problems,
corresponding to equations (\ref{equation-IVP-D1})
and (\ref{equation-IVP-D2}). Given $a\in J$ and $\alpha\in J$,
there are unique $\lambda$ and $\mu$, repectively,
 such that the solutions
$\Phi_+(a,\alpha)$ and
$\Phi_-(a,\alpha)$, respectively,
to these equations conjugate $f$ to $g$ on $J$ and
have $C^1$ extensions to $(c,d]$ and $[c,d)$,
respectively.  We may then formulate a solution to
the conjugacy problem, as follows:

\begin{theorem}\label{theorem-compact}
Let $I$ be a compact interval and
let $f,g\in\Diffeo(I)$, both fixed-point-free on $J$, both in $S_-$.
Then the following conditions are equivalent:\\
(1) $f$ is conjugate to $g$ in $\Diffeo(I)$;
\\
(2) The product $H(x,\xi)$ converges for some (and hence for
all) $x>0$ and $\xi>0$, and there exists some $a>0$
and $\alpha>0$ such that the solution $\Phi_+(a,\alpha)$
extends $C^\infty$ to both ends of $I$;
\\
(3) There exist $a>0$ and $\alpha>0$ such that $H(a,\alpha)$
converges, and $\Phi_+(a,\alpha)=\Phi_-(a,\alpha)$
extends in $\Diffeo(I)$.
\qed
\end{theorem}

Details are in Section \ref{section-compact}.

\subsection{Compact $I$: Functional Moduli}
For some problems of classification, a solution is available in
terms of a finite-dimensional space of invariants, or
\lq\lq moduli". But if the class structure is very rich, this may
not be possible, and only infinite-dimensional spaces of moduli are
natural. This is the genesis of the idea of functional modulii
(cf. [V]).

In special cases, the conjugacy problem on a compact interval
can be reduced to condition (T) at both ends, plus identity
of a suitable modulus (a conjugacy invariant that is a
diffeomorphism on some interval).  See Robbins \cite{RO}, Afraimovitch
Liu and Young \cite{ALY}, and Young \cite{Y}.
All these results are subsumed in an unpublished
lemma of Mather \cite{M2},
subsequently and independently found by Young,
which covers the case in which the germs of $f$ at both ends
of the interval are the exponentials of smooth vector fields, and
for which the modulus is a double coset $RkR$ of the rotation group
$R=$SO$(2,\R)$
in the group $\Diffeo^+(\S^1)$ of circle diffeomorphisms, and the conjugacy
class of $f$ is determined by the smooth conjugacy classes
of the two vectorfield germs and the modulus.

See Subsection \ref{subsection-modulus} for more detail on moduli.

\subsection{Compact $I$: Shape}\label{subsection-shape}
Obviously, it is rather unlikely that two given maps
$f$ and $g$ will be smoothly conjugate on $I$,
even assuming they satisfy the sign condition and Condition (P).

The conditions of Theorem \ref{theorem-compact}
 are necessary and sufficient, but are
tedious to check.

It is worth noting a necessary
condition (the \lq\lq shape" condition)
that is easier to check in the compact case. This
will often suffice to
show two maps are not conjugate.

First we
define
$$ F_{a}(x) = H(f,f;x,a)=
\prod_{n=-\infty}^\infty
\frac{f^\prime(f^{\circ n}(x))}%
{f^\prime(f^{\circ n}(a))},$$
and
$$ G_{\alpha}(\xi) = H(g,g;\xi,\alpha) =
\prod_{n=-\infty}^\infty
\frac{g^\prime(g^{\circ n}(\xi))}%
{g^\prime(g^{\circ n}(\alpha))}$$
whenever $x,\xi,a,\alpha\in J$.
Note that
\begin{equation}\label{equation-FHHG}
 H(x,\xi) \cdot G_\alpha(\xi) = F_a(x)\cdot H(a,\alpha),
\end{equation}
whenever all the terms make sense.

\begin{proposition}\label{proposition-shape}
Suppose $f,g,h\in \Diffeo(I)$, $f$ is fixed-point-free on $J$, and
$f=g^h$. Then $H(x, h(x) )$ is constant on $J$.  Thus, given any $a,
\alpha\in J$, there is some $\kappa>0$ such that
$$ F_a(x) = \kappa G_\alpha(h(x)), \ \forall x\in J. $$
\end{proposition}
\begin{proof}
Suppose $I$ is compact, $f\in S_-$,  and
$f=g^h$ in $\Diffeo(I)$.
Applying the results about half-open intervals
to both $[c,d)$ and $(c,d]$, we see
that the product $H_1(x, h(x))$ converges to $h^\prime(x)/h^\prime(d)$
for each $x\in J$, and
$H_2(x, h(x))$ converges to $h^\prime(x)/h^\prime(c)$
for each $x\in J$.
Thus the two-sided product
$$H(x, h(x)) = H_1/ H_2 =
\prod_{n=-\infty}^\infty
\frac{ f^\prime(f^{\circ n}(x))}
{g^\prime(g^{\circ n}(h(x)))} $$
is independent of $x\in J$, and equals the ratio
$h^\prime(c)/h^\prime(d)$ of the derivatives of
the conjugating map at the ends.

This immediately tells us
that $H(x,h(x))$ is constant. The rest then follows from
equation (\ref{equation-FHHG}).
\end{proof}

This means that the graphs of each $F_a$ and of each $G_\alpha$ have
the same \lq\lq shape".   If they are not monotone, then the
relative diffeomorphism class of the critical set and the pattern of
maxima and minima must be the same for both functions. The pattern
for $F_a$ is determined by the pattern on the segment
$I_a=[a,f(a)]$, because it repeats itself on successive images of
$I_a$ under $f$. Similarly, the pattern for $G_\alpha$ is determined
by the pattern on $[\alpha, g(\alpha)]$. Apart from this
quasiperiodic feature, the patterns may be pretty complicated.

Note that if the condition of the proposition fails,
then this can be determined by a computation.

\subsection{Compact $I$: Flowability}\label{subsection-flowability}
We note applications to existence of a smooth
flow on a compact interval $I=[c,d]$,
for which $f$ is the time-1 step.

Applying Theorem \ref{theorem-compact} to the case
$g=f$, we see that the centraliser
$C_f$ is the intersection of two at-most-one-parameter
groups, containing the (nontrivial, discrete) group
of all iterates of $f$.
(It may well be that only the compositional powers of
$f$ belong to $C_f$.  )

We deduce a method for
deciding whether or not $f$ is the time 1 map
of a flow.

\begin{proposition}\label{proposition-flow}
A diffeomorphism $f\in\Diffeo^+(I)$
is flowable if and only if the
centralisers of $f$ in $\Diffeo(\{c\}\cup J)$
and $\Diffeo(J\cup\{d\})$ are both connected,
and coincide
(when restricted to $J$).
\qed
\end{proposition}

Applying the shape result,
Proposition \ref{prop-flow-shape}, we
identify a special case in which a necessary condition
for flowability may be checked by plotting a graph.

\begin{proposition}\label{prop-flow-shape}
Suppose that $f\in\Diffeo(I)$, is fixed-point-free on $J$
and $f$ is flowable. Then
for each $a\in J$, $F_a$
is either strictly monotone on $J$, or constant on $J$.
\end{proposition}
\Proof Suppose that $F_a$ is neither strictly monotone
on $J$ nor constant on $J$.
Each conjugacy if $f$ to itself must permute the maximal open
intervals of strict monotonicity of $F_a$. Since $F_a$ is smooth and not
strictly monotone or constant,
there exist at least two such intervals, and since the
pattern repeats, there are in fact infinitely many. But the number
is countable, since they are pairwise disjoint open sets, and
conjugacy must permute the countable set of endpoints of these
intervals of monotonicity, and is determined uniquely by the image
of one endpoint. Hence the centralizer of $f$ is a countable group,
so $f$ cannot be flowable.
\qed

We can do better when the graph of $f$ is tangent to the diagonal at
the ends of $I$:

\begin{corollary} Suppose $f\in\Diffeo(I)$ is fixed-point free on
$J$ and is flowable. Then the following are equivalent:

\begin{enumerate}
    \item $f'(c)=f'(d)$;
    \item $f'(c)=f'(d)=1$;
    \item $F_a$ is constant on $J$, for each (or any one) $a\in J$.
\end{enumerate}
\end{corollary}
\Proof The implication $(1)\implies(2)$ follows from the fact that
$1$ is always trapped between $f'(c)$ and $f'(d)$.

Next, note that we have the formula
\begin{equation}\label{equation-Fafx}
F_a(f(x))=F_a(x)\frac{f'(d)}{f'(c)},
\end{equation}
whenever $a,x\in J$.

Suppose $(2)$ holds. Fix $a\in J$. The formula (\ref{equation-Fafx})
implies that
$F_a(f(x))=F_a(x)$ for all $x\in J$. Since $f$ is flowable,
Proposition \ref{prop-flow-shape} tells us that $F_a$
is constant on each interval $[f(x),x]$. But for any fixed
$x_0=x$, the iterates $x_n=f^{\circ n}(x)$ converge
monotonically to one end of $J$ as $n\uparrow+\infty$,
and monotonically to the other end as $n\downarrow-\infty$,
hence the intervals $[x_{n+1},x_n]$ pave $J$, and,
since $F_a$ is constant on each, it is constant on the
whole interval $J$.  Thus $(2)\implies(3)$.

Finally, suppose $(3)$ holds. Then equation (\ref{equation-Fafx}),
applied to any $x\in J$,
yields $f'(d)=f'(c)$, since $F_a(x)$ never vanishes.
\qed

We note that these results depend only on the assumption that $f\in
C^2(I)$.

\subsection{Conjugacy in $\Diffeo^+(I)$}\label{subsection-Diffeo+}
Now we move on to the general orientation-preserving case on an arbitrary
interval $I\subset \R$.

Each interval is diffeomorphic to one of the closed intervals $\R$,
$[0,+\infty)$ or $[-1,1]$, so there is no loss in generality
in supposing that $I$ is
a closed interval. (If $I$ is not closed, fix some diffeomorphism
$h$ of $I$ onto a closed interval. Then $f$ and $g$ are conjugate
in $\Diffeo^+(I)$ if and only if $^hf$ and $^hg$ are conjugate
in $\Diffeo^+(h(I))$.)

We can reduce the problem to the conjugacy
problem in
\\ $\Diffeo^+_B=\Diffeo^+_B(I)$, with $B=\bdy E$,
and
$E\subset I$ closed:
\begin{proposition}\label{proposition-reduction-2}
Suppose $I$ is a closed interval.
Let $f,g\in\Diffeo^+$. Then there exists $h\in\Diffeo^+$
such that $f=g^h$ if and only if
there exists $h_1\in\Diffeo^+$ such that the following three conditions
hold:
\\
(1) $h_1(\fix(f))=\fix(g)$;
\\
(2) letting $f_1=f^{h_1^{\circ-1}}$, and $E=\fix(g)$, we have
$$ \sign(f_1(x)-x) = \sign(g(x)-x),\ \forall x\in I\sim E;$$
(3) there exists $h_2\in\Diffeo^+_{\bdy E}$ such that
$f_1 = g^{h_2}$.
\end{proposition}

\Proof \lq\lq Only if": Suppose
there exists $h\in\Diffeo^+$
such that $f=g^h$.

Taking $h_1=h$, we have condition (1).
Also, in that case $f_1=g$, so condition (2) holds. Taking
$h_2=Id$, we have condition (3).

\lq\lq if": Suppose there exist $h_1$ and $h_2$
satisfying conditions (1), (2) and (3).
Then $h=h_2\circ h_1$ has $f=g^h$.
\qed

As already remarked, the existence of an $h$ with condition
(1) is not amenable to algorithmic checking, so we
shall just treat it as a primitive condition.

Given the existence of such $h$, there may exist many. Condition
(2) cuts down the collection of eligible $h$. One then has to
check condition (3) for each eligible $h$.
In this sense, we have reduced the conjugacy problem in
$\Diffeo^+$ to the problem of characterising conjugacy
in $\Diffeo^+_{\bdy E}$, for two elements of
$\Diffeo^+_{E}$
(-- see Example \ref{example-2.2}).
It is worth remarking that the condition that $T_af$ and $T_{h(a)}g$
be conjugate Taylor series reduces the collection of eligible $h$
considerably.
See further remarks in Subsection
\ref{examples-subsection-Diffeo^+}.

Obviously, the reduction achieved here is not deep.  However it is useful.
If we replace $\bdy E$ by $E$ in condition (3), then the restated proposition
remains true, but is less useful.  To explain this point, consider
this example.

\begin{example}\label{example-2.2}
\end{example}
Take $I=\R$ and set
$$ g(x) = \left\{
\begin{array}{rcl}
x &,& x\le0,\\
x+e^{-1/x} &,& x>0,
\end{array}
\right.
$$
and $f(x)=\half(1+g(2x-1))$.  Then
$x\mapsto 2x-1$ conjugates $f$ to $g$
in $\Diffeo^+$.

If we are handed two functions $f$ and $g$ in $\Diffeo^+$, and
asked to determine whether or not they
are conjugate, then we would begin by comparing the
pairs $(\R, \bdy(\fix(f)))$ and $(\R,\bdy(\fix(g)))$, to see whether
they lie in the same diffeomorphism class.  In the present example,
the pairs are $(\R,\{\half\})$ and $(\R,\{0\})$, and (of course)
they do.
The next reasonable
step would be to take {\em any} $h_1\in\Diffeo^+$ that maps
$\half$ to $0$, and replace $f$ by
$f_1=f^{h_1^{\circ-1}}$, as in condition (2) of the proposition,
and proceed to compare $f_1$ and $g$, which now have the same
fixed point set, namely $E=(-\infty,0]$.  Let's say we chose
(slightly perversely),
$ h_1(x)=4x-2$.  Then we would have $f_1(x)=g(x)=x$ for $x\le0$,
and $f_1(x) = 2g(x/2)$ for $x>0$.  We would then proceed to
check $f_1$ and $g$ are conjugate, and it would be reasonable
to seek a conjugacy on the lines of (3), that fixes $0$.
We could then use the differential equation, as in Section
\ref{section-half-open} to \lq\lq discover"
one of the conjugacies that exist in $\Diffeo([0,+\infty))$
(i.e. one of the elements of
the coset of $C_f$ to which $x\mapsto 2x$ belongs).  Now each such conjugacy
is $C^\infty$ down  to $0$, and has derivative $2$ at $0$.
Extending it {\em in any way whatsoever} to a diffeomorphism
of $(-\infty,0]$ gives a global conjugacy $h_2$
from $f_1$ to $g$, because
both maps are the identity map on the negative axis, and are conjugated
by anything.  Each $h_2$ found in this way fixes $\{0\}$  but,
and this is the point, {\em there is no
conjugacy of $g$ to $f_1$ that fixes all points of $E$.
 In order to attack the problem in this way
it is essential to retain the flexibility to move points inside
$\fix(g)$}.  Otherwise, this approach goes nowhere.

\subsection{Conjugacy in $\Diffeo^+_B$}\label{subsection-Diffeo+B}
{\bf Throughout this subsection, $I$ will be a fixed
closed interval (bounded or not), $E$ will be a fixed closed
nonempty subset of $I$, containing all endpoints of $I$,
and $B$ will be the boundary of $E$.}

\medskip
We consider $f,g\in \Diffeo^+(I)$,
(Recall that  in view of Proposition \ref{proposition-fpfree}
we need not consider the special case $E=\emptyset$.)

As before, we suppress the explicit $(I)$ in
$\Diffeo^+(I)$, $\Diffeo^+_E(I)$, $\Diffeo^+_B(I)$, etc..
$I\sim E$ is a countable union of open
intervals.
The following is trivial:
\begin{proposition}
Let $f,g\in\Diffeo^+$, fixing precisely $E$. Then
$f$ is conjugate to $g$ in $\Diffeo^+_B$
if and only if
there is a  global function
$h\in\Diffeo^+$ such that, for each connected
component $J$ of $I\sim E$, the restriction to
each clos$(J)$ belongs to $\Conj(f,g;J)$.
\end{proposition}
\qed
So necessary conditions for the conjugacy
are:

\medskip
(1) If $J$ is an unbounded component of $I\sim E$,
then the restrictions of $f$
and $g$ to the closure of $J$ satisfy
the conditions of
Corollary \ref{corollary-half-open-}
or Corollary \ref{corollary-half-open+},
as appropriate (i.e. depending on
whether or not $f$ iterates points towards
or away from the (finite) end of $J$);
\\(2)
If $J$ is an bounded component of $I\sim E$,
then the restrictions of $f$
and $g$ to the closure of $J$ satisfy
the conditions of Theorem \ref{theorem-compact};
\\(3)
If an endpoint $p=c(J)$ or $d(J)$ is not isolated in $B$,
then some conjugating $h\in\Diffeo^+_B$ must have
$h-x$ flat at $p$.

Condition (3)
actually implies that
{\em all} elements
of $\Conj(f,g;J)$ must be flat at $p\,$:

\begin{lemma}\label{lemma-Diffeo+B-flat}
Suppose $f$ is conjugate to $g$ in $\Diffeo^+_B$.
Then whenever an end $c(J)$ or
$d(J)$ of a component $J$ is an accumulation point of $B$, it
follows that each $\phi\in \Diffeo^+(\clos(J))$ that conjugates $f$ to
$g$ on $I$ has $\phi(x)-x$ flat there.
\end{lemma}

The proof is in Section \ref{section-Diffeo+B}.

So we formulate this version of condition (3):

\medskip
\noindent {\bf Condition (F):}
{\sl If $J$ is a connected component of $I\sim E$ and
an endpoint $p=c(J)$ or $d(J)$ is not isolated in $B$,
then each conjugating $h\in\Diffeo^+_B$ must have
$h-x$ flat at $p$.}

\medskip
But the conditions (1)-(3) will not
always guarantee the existence of
a global conjugation,
even when $E$ is finite, or,
more generally, discrete.
The family $\Conj(J)=\Conj(f,g;J)$ is mapped by $h\mapsto T_ph$
to a set $M(p,J)=M(f,g;p,J)$ of Taylor series, whenever $p$
is a finite end of $J$. This set $M(p,J)$ is a coset of
a subgroup of the invertible Taylor series, and is,
generically, discrete.  Whenever two intervals
$J$ and $J'$ have a common endpoint $p$,
we are snookered unless $M(p,J)$ and
$M(p,J')$ intersect.  This gives us a
necessary condition:

\medskip
\noindent
{\bf Condition (M$_1$):} {\sl If $p$ is an isolated point of $B$,
and $J$ and $J'$ are the components of $\R\sim E$
to the left and right of $p$, then
$M(p,J)\cap M(p,J')\not=\emptyset$.}

\medskip
In case $p$ is a hyperbolic fixed point, Condition (M$_1$) is
equivalent to the simpler condition that the multiplier cosets
$\{h'(p):h\in\Conj(f,g,J)\}$ intersect. (These are cosets of a
subgroup of the multiplicative group $(0,\infty)$.) For suppose we
take conjugating maps $h$ and $k$ on clos$J$ and clos$J'$,
respectively, with the same multiplier at $p$. Then $T_ph$
conjugates $T_pf$ to $T_pg$, and so does $T_pk$, so
$(T_ph)\circ(T_pk)^{\circ-1}$ commutes with $T_pf$ and has
multiplier $1$.
Now a hyperbolic series is conjugate
to its linear part, and the centraliser of a linear series
is the set of all linear series, so an element of
the centraliser
of a hyperbolic series is determined uniquely by its multiplier.
Hence
$(T_ph)\circ(T_pk)^{\circ-1}$
equals $X$.

A similar argument shows that if
$T_pf = X$ mod$X^p$, but $T_pf\not=X$ mod$X^{p+1}$,
then Condition (M$_1$) simplifies to the condition that
the quotients mod $X^{p+1}$ intersect:
$$\left(M(p,J)\mod X^{p+1}\right)
\cap \left(M(p,J')\mod X^{p+1}\right)\not=\emptyset.$$

\medskip
\begin{example}
\end{example}
Take $f$ and $g$, fixing only $0$, with $f-x$ and
$g-x$ flat at $0$, such that $f$ is conjugated
to $g$ on $J_1=[0,+\infty)$ by $x\mapsto 2x$,
and $f=g$ on $J_2=(-\infty,0]$.
Then $\Conj(f,g,J_1)$
is nonempty, but only has maps with derivative $2$ at $0$,
whereas $\Conj(f,g,J_2)$ has only maps with derivative $1$
at $0$, so condition (M$_1$) fails.

\bigskip
Assuming Condition (M$_1$), we have a further problem
if there is a chain of successive isolated points in $E$.
We then have a chain $J_1,\ldots J_k$ of
successive components of $I\sim E$. To find a conjugation,
we must patch together elements of the $\Conj(f,g;J_i)$
to make a single smooth conjugation on
the closure of the union of the $J_i$. If we begin with one
element of $\Conj(f,g;J_1)$, and work along, trying to match its
Taylor series at each common endpoint,
then it becomes increasingly improbable that we will succeed.
If there is any chain of intervals for which it cannot be
done, then there is no global conjugation.

The key to further progress is to focus on
$B'$, the set of accumulation points of $B$.
The connected components of $I\sim B'$ include
the connected components $J$ of the interior
of $E$, and on these $J$ {\em every} diffeomorphism
conjugates $f$ to $g$. This makes it reasonable to define
$\Conj(f,g;J)=\Diffeo^+(\clos J)$ for such $J$.

We formulate a stronger version of condition(M$_1$):

\medskip
\noindent
{\bf Condition M$_2$:}
{\sl Given any connected component $L$ of $\R\sim B'$, there exists
a function $h\in\Diffeo^+(L)$ whose restriction to
each component $J$ of $L\sim B$ belongs to $\Conj(f,g;J)$.}

\medskip
Each set $L\cap B$ is empty, or finite, or forms a sequence
tending to one end of $L$, or a two-sided sequence
accumulating at both ends of $L$.

If the condition fails, then it can be disproved
by starting in any one $J$, and working left and right,
cutting down the set of eligible conjugations,
until at some stage it is found that the set is
empty.

However, it may be that all the functions that thread together
conjugations on the various $J$ wiggle too much to
extend smoothly to the accumulation points at the ends.
So we formulate:

\medskip
\noindent
{\bf Condition (M):}
{\sl Given any connected component $L$ of $I \sim B'$, there exists
a function $h\in\Diffeo^+(\clos L)$ whose restriction to
each component $J$ of $L\sim B$ belongs to $\Conj(f,g;J)$.}

Note that this implies conditions (1) and (2).

It may seem that we are heading into
a jungle as complex as that involved in the
order-equivalence problem, and that
higher derived sets are about to appear.
But the situation is not so bad. We do not have to look at $B''$:

\begin{theorem}\label{theorem-Diffeo+B}
Let $I$ be a closed interval. Let
$f,g\in\Diffeo^+(I)$ both fix precisely $E\subset I$. Then
$f$ is conjugate to $g$ in $\Diffeo^+(I)$ (or, equivalently,
in $\Diffeo^+_B(I)$)
if and only if both conditions (M) and (F) hold.
\end{theorem}

The proof is in Section \ref{section-Diffeo+B}.

\subsection{Reducing from $\Diffeo(I)$ to $\Diffeo^+(I)$}\label{subsection-Diffeo}
In Section \ref{section-Diffeo}
we discuss the reduction of the conjugacy
problem in  the full diffeomorphism group
to the conjugacy problem in the subgroup
of direction-preserving maps.

We close with some examples in Section \ref{section-examples}.

\section{The Fixed-Point-Free Case}\label{section-fpfree}
In this section our main purpose is to prove
Proposition \ref{proposition-fpfree}:
Some of the elements of the proof will be useful later,
for other purposes.

So suppose $f$ and $g$ are fixed-point-free
elements  of $\Diffeo^+(\R)$, and both move all points up, or both
move all points down. Then  we have to show that $f$
and $g$ are conjugate in
$\Diffeo^+(\R)$.

The proof depends on a well-known result due to \'E. Borel(cf. \cite{M}, or \cite{R},
Chapter 19]):
\begin{theorem} For each $a\in\R$, each formal
power series is the power series at $a$ of some smooth function.
\end{theorem}
\qed
\begin{corollary}\label{corollary-3.2}
Given a point $a\in\R$, any value $\lambda\in\R$, and a power series
$P=a_1X+\cdots$ with $a_1>0$, there exists $f\in\Diffeo^+$ with
truncated Taylor series $T_af=P$, and with $f(a)=\lambda$.
\end{corollary}
\Proof
First, pick a smooth function $h_1$
with Taylor series at $a$ equal to the
term-by-term derivative $P^\prime$ of $P$.
Then $h_1$ will be positive near $a$, so by
modifying it off a neighbourhood of $a$
one may construct an everywhere-positive
smooth function $h_2$ with the Taylor series
$P^\prime$ at $a$.  Now take
$$ f(x) = \lambda+\int_a^x h_2(t)\,dt, \ \forall x\in\R.$$
\qed
\begin{corollary}
Given real numbers $a<b$ and formal series
$$ P = a_1X+\cdots,\qquad Q = b_1X+\cdots
$$
with $a_1>0$, $b_1>0$, there exists
$f\in\Diffeo$ with
$$ f(a)=a,\ f(b)=b,\ T_af=P, \ \hbox{ and }T_bf=Q.$$
\end{corollary}
\Proof  Applying the previous corollary twice, choose
diffeomorphisms $r$, and $s$ such that
$$\begin{array}{rcl}
r(a)&=a,&\quad  T_ar=P,\\
s(b)&=b,&\quad  T_bs=Q.
\end{array}
$$
Since $r^\prime(a)>0$ and $s^\prime(b)>0$, we may choose $\eta>0$
such that $r$ maps $[a,a+2\eta]$ diffeomorphically onto
$[a,r(a+2\eta)]$, $s$ maps $[b-2\eta,b]$ diffeomorphically onto
$[s(b-2\eta),b]$, and
$$\max\{a+2\eta,r(a+2\eta)\}<\min\{b-2\eta,s(b-2\eta)\}.$$

Choose a monotonically nonincreasing
 smooth function $t$ that is identically $1$
on $(-\infty,a+\eta]$ and is identically $0$ on
$[a+2\eta,+\infty)$.

Choose a monotonically nondecreasing
 smooth function $u$ that is identically $0$
on $(-\infty,b-2\eta]$ is identically $1$ on
$[b-\eta,+\infty)$.

Choose another smooth function $v$ that is nonnegative,
is {\em not} identically zero, but is zero off
$[a+\eta,b-\eta]$.

For each $\lambda>0$, let
$$
h_\lambda(x) = t(x)\cdot r^\prime(x) + u(x)\cdot s^\prime(x)
+\lambda\cdot v(x),\ \forall x\in\R.$$

Then $h_\lambda(x)>0$ for all $x\in\R$,
$h_\lambda(x)=r^\prime(x)$ whenever $x<a+\eta$,
$h_\lambda(x)=s^\prime(x)$ whenever $x>b-\eta$.

Define
$$ f_\lambda(x) = a + \int_a^x h_\lambda(z)\,dz,\quad\forall x\in\R.$$

Then $f_\lambda$ is a diffeomorphism, and fixes $a$.
Also $f_\lambda$ has truncated Taylor series
$P$ at $a$ and $Q$ at $b$. To finish, we just need to
pick $\lambda>0$ so that $f_\lambda(b)=b$.
This can be done, because it amounts to solving
$$ \lambda\int_a^bv(z)\,dz = b-a - \int_a^{a+2\eta}t(z)r^\prime(z)\,dz
-\int_{b-2\eta}^bu(z)s^\prime(z)\,dz$$
and the right-hand side is positive, since
it exceeds
$$ b-a - \int_a^{a+2\eta}r^\prime(z)\,dz
-\int_{b-2\eta}^bs^\prime(z)\,dz
= s(b-2\eta)-r(a+2\eta)>0.$$
\qed

Now we can give the proof of Proposition \ref{proposition-fpfree}:

\Proof
Since $\R$ is diffeomorphic to each nonempty open interval
$I\subset\R$, it suffices to prove
Proposition \ref{proposition-fpfree} for the case $I=\R$.

By means of a preliminary conjugation with a linear map,
we may assume that $g(0)=f(0)$. Let $b=f(0)$.

Pick a smooth increasing map $\varphi$ of the interval
$ [0,b]$ onto itself such that
$$ (T_b\varphi)\circ (T_0 f)=(T_0 g)\circ (T_0\varphi).$$
(For instance, one could take $T_0\varphi\equiv X$, and let the
above equation define $T_b\varphi$; the existence of a
$\varphi$ matching these Taylor series follows from the last Corollary.)

The conjugacy equation then forces a unique extension of $\varphi$
to an element of $\Diffeo^+(\R)$.
\qed

We see that there are a great many conjugacies between two given
conjugate diffeomorphisms. In particular, the centraliser of
a fixed-point-free diffeomorphism is very large, and is not
abelian. We shall see below that the presence of even a single
fixed point produces a drastic reduction in the size
of the centraliser.  It becomes at most a one-parameter abelian group.

\section{$[0,+\infty)$: Necessary conditions.}\label{section-half-open}
\subsection{Proof of Theorem \ref{notone}}\label{SLT}
\Proof
Obviously, (2) implies (1).

{(1) implies (3):} Let $f'(0)=g'(0)=a$. Observe that, as $f\in
S_-$, $0<a<1$. By Sternberg there exist $\phi, \psi\in
\Diffeo^+([0,\infty))$, such that $\phi^{-1}\circ f\circ \phi(x)=a
x=\psi^{-1}\circ g\circ \psi(x)$. Now, if $\tau_a(x)=a x$, the
functions of the sequence can be presented as $h_n=\psi\circ
\tau_{a}^{-n}\circ \psi^{-1}\circ \tau_\lambda \circ \phi \circ
\tau_a^n \circ \phi^{-1}$. As both $\phi$ and $\psi$ are
diffeomorphisms, $\phi(x)=\phi'(0)x+O(x^2)$, and $\psi^{-1}(x)=
(\psi^{-1})'(0)x+ O(x^2)$. This means that $\phi \circ \tau_a^n
\circ \phi^{-1}(x)=\phi'(0) a^n \phi^{-1}(x)+ O((a^{n}
\phi^{-1}(x))^2)$, and $\psi^{-1}\circ \tau_\lambda \circ \phi \circ
\tau_a^n \circ \phi^{-1}=(\psi^{-1})'(0)\lambda \phi'(0) a^n
\phi^{-1}(x)+ O((a^{n} \phi^{-1}(x))^2)$. Placing this expression in
the formula for $h_n$ we get $h_n(x)=\psi((\psi^{-1})'(0)\lambda
\phi'(0) \phi^{-1}(x)+O(a^n\phi^{-1}(x)))$. For a fixed $x$ we see
that $h_n(x)\rightarrow \psi((\psi^{-1})'(0)\lambda \phi'(0)
\phi^{-1}(x))$, when $n\rightarrow \infty$, which is a
diffeomorphism.

Obviously, (3) implies (4), and (4) implies (5).

{(5) implies (2):} We see that $g\circ h(x)=g(\lim(h_n(x)))=
\lim g(h_n(x))=\lim h_{n-1} (f(x))=h\circ f(x)$, and $h\in
\Diffeo^+([0,\infty))$.
\qed

The nice thing about this is that (3)-(5) give us
a constructions for conjugacy maps. Ahern and Rosay call
the construction (4) \lq\lq the basic trick".

\subsection{Proof of Theorem \ref{theorem-2.4}}\label{subsection-takens}
Any two elements of
$S_-$ that
agree on a neighbourhood of $0$
are smoothly conjugate. This is easy to see: one just uses the
conjugacy equation to extend the trivial conjugation given
by the identity function near $0$ to a smooth conjugacy on the
whole of $[0,+\infty)$.

So to prove Taken's Conjugacy Theorem
we just have to show that if $f\in\Diffeo([0,+\infty))$ and
the series $T_0f$ is conjugate to
$X-X^{p+1}+\alpha X^{2p+1}$, then
$f$ is smoothly conjugate to a diffeomorphism
that coincides with
$g=x-x^{p+1}+\alpha x^{2p+1}$
for small enough $x>0$.

However, the constructive method of
the proof of Theorem \ref{notone} cannot be used
directly in this case, in order to find a conjugacy.

\begin{example}
\end{example}

Consider
$$ f(x)= x+x^2,\ \forall x\ge0. $$
This function belongs to $S_+$ and
is conjugate to
$$ g(x) = x + 2x^2 = \half f(2x),$$
and hence the conjugating map $h:x\mapsto \frac12 x$ conjugates
$f^{\circ-1}\in S_-$ to $g^{\circ-1}$. But the conjugation functions
in general cannot be recovered as in condition 3 of the Theorem
\ref{notone}. For instance, taking $\lambda=2$, we find
$$ g^{\circ n}(2f^{\circ-n}(x)) \rightarrow \infty,$$
and does not converge at all. One can see this by observing that
$$
g^{\circ n}(2f^{\circ-n}(x))=\frac12 f^{\circ n}(4f^{\circ-n}(x)),$$
due to the conjugation. On the other hand given a small number $x$,
the first number $N$ for which $f^{\circ N}(x)\geq 4x$ is at least
$\frac{3x}{(4x)^2}=\frac3{16x}$, and as the functions are monotone,
we see that $f^{\circ n}(4f^{\circ-n}(x))\geq f^{\circ (N_n)}(x)$,
where $N_n=\lfloor\frac3{16 f^{\circ -n}(x)}\rfloor\rightarrow
\infty$ when $n\rightarrow \infty$. Thus $ g^{\circ
n}(2f^{\circ-n}(x))\geq \frac12 f^{\circ (N_n)}(x)\rightarrow
\infty$.

\medskip
Takens proceeds in two steps.
\\
(1) He shows that $f$ is conjugate to a diffeomorphism
$f_1(x)= x-x^{p+1}+\alpha x^{2p+1}+g_1(x)$, where
$g_1$ is a $C^\infty$ function flat at $0$.
We can see this at once from Corollary \ref{corollary-3.2}:
The Taylor series of $f$ is conjugate to $X-X^{p+1}+\alpha X^{2p+1}$,
so choose an invertible series $H$ that has
$$ H^{\circ-1}T_0fH = X-X^{p+1}+\alpha X^{2p+1}.$$
Then choose $h\in\Diffeo(\R)$i, fixing $0$,
with $T_0h=H$, and let $f_1 = f^h$.
\\
(2) He defines $\Psi:\R^2\to\R^2$ by
$$ \Psi(x,t) = (x-x^{p+1}+\alpha x^{2p+1}+tg_1(x),t).$$
and shows there is a vectorfield
$$ \tilde Z = Z(x,t)\frac{\partial}{\partial x}
+\frac{\partial}{\partial t} $$
with $Z$ flat at all points where $x=0$, and
$\Psi_*(\tilde Z) =\tilde Z$ on a neighbourhood
of the segment $\{(x,t):x=0,0\le t\le1\}$.
The proof of the existence of such a vectorfield
(which is a fixed-point theorem) requires a
substantial argument (\cite[pp. 177-189]{T})
and we omit the details. Once he has it,
the conjugacy of $f_1$ to $g$ is
obtained by taking $\phi(x)$ so that,
for small $x>0$, the points $(x,1)$
and $(\phi(x),0)$ lie on the same integral curve
of $\tilde Z$.
\qed

\medskip
Ahern and Rosay \cite[pp. 549-51]{AR} give another
proof of his theorem.  They show that, in fact, if
condition (T) holds, then the
\lq\lq basic trick" construction of Theorem \ref{notone}, condition (4)
may be used, with caution, on a neighbourhood of
$0$, to get a conjugacy started.

\section{$[0,+\infty)$: The Product Condition}
\label{section-necessary}
Fix arbitrary $f,g\in S_-$.

\begin{lemma}\label{lemma-basic-product} Suppose $f$ and $g$ are conjugate in
$\Diffeo([0,\infty))$. Then for any $x>0$ there exists $\xi>0$
such that
the product (\ref{equation-basic-product})
converges.
\end{lemma}
\Proof Pick $h\in \Diffeo([0,\infty))$ with $f=g^h$, and
set $\xi=h(x)$. We observe that $h\circ f^{\circ n}=g^{\circ n}\circ
h$, hence equating derivatives we get
$$
h'(f^{\circ n}(x))\frac{df^{\circ n}}{dx}(x)=\frac{dg^{\circ
n}}{d\xi}(\xi)h'(x),
$$hence
$$\prod\limits_{j=0}^{n-1}\frac{f'(f^{\circ j}(x))}{g'(g^{\circ j}(\xi))}=
\frac{h'(x)}{h'(f^{\circ n}(x))},$$
so the product converges to the
limit $h'(x)/h'(0)$.
\qed

The correspondence between $x$ and $\xi$, referred to in the
lemma
is not essential, for we have the following, which is due to
Kopell \cite{K}.  (We give the proof for convenience.)

\begin{lemma}\label{lemma-product-existence}
 Let $x,y\in [0,\infty)$ and denote $x_n=f^{\circ
n}(x)$, $y_n=f^{\circ n}(y)$. Then the infinite product
\begin{equation}\label{product-just-f}
\prod\limits_{n=0}^{\infty}\frac{f'(x_n)}{f'(y_n)}
\end{equation}
converges.
\end{lemma}

\Proof
First, assume that $y_0$ is between $x_1$ and $x_0$. The convergence
of the product is equivalent  to the convergence of the series of
logarithms $\sum\limits_{n=0}^{\infty}\ln(\frac{f'(x_n)}{f'(y_n)})$,
which in turn is equivalent to that of $\sum\limits_{n=0}^{\infty}
\left(1-\frac{f'(x_n)}{f'(y_n)}\right)$. Now
$$\left|\frac{f'(x_n)-f'(y_n)}{f'(y_n)}\right|\leq
\left(\frac{\sup|f''|}{\inf|f'|}\right)\cdot|x_n-y_n|$$
(where the
$\sup$ and $\inf$ are taken on $[0,x]$; note that the $\inf$
is positive since $f$ is a diffeomorphism), and so the convergence
follows from $\sum\limits_{n=0}^{\infty}|x_n-y_n|\leq |x_0|$, which
holds because the intervals from $x_n$ to $y_n$ are
pairwise-disjoint subintervals of that from $0$ to $x_0$.

So the result holds when $y_0$ is between $x_1$ and $x_0$.

For general $y$, choose $k\in\Z$ such that
$y_0$ is between $x_k$ and $x_{k+1}$.
Then
$$
\prod\limits_{n=0}^{m}\frac{f'(x_n)}{f'(y_n)}
=
\frac{f'(x_0)\cdots f'(x_{k-1})}{f'(y_{m-k+1})
\cdots f'(y_m)}\cdot
\prod\limits_{n=0}^{m-k}\frac{f'(x_{n+k})}{f'(y_n)},
$$
so the series
$$ \prod\limits_{n=0}^{\infty}\frac{f'(x_n)}{f'(y_n)}
\hbox{ and }
\prod\limits_{n=0}^{\infty}\frac{f'(x_{n+k})}{f'(y_n)}
$$
converge or diverge together, with
$$
\prod\limits_{n=0}^{\infty}\frac{f'(x_n)}{f'(y_n)}
=
\frac{f'(x_0)\cdots f'(x_{k-1})}{f'(0)^k}\cdot
\prod\limits_{n=0}^{\infty}\frac{f'(x_{n+k})}{f'(y_n)}.
$$
Thus we obtain the general result by replacing
$x_n$ by $x_{n+k}$.
\qed

\begin{corollary}\label{corollary-flat-ff-product}
(1) In case $x_1<y<x_0$, and
$T_0f=X+bX^{p+1}+\ldots$ for some $p\in\N$,
the product (\ref{product-just-f})
is $1+\Oh(x^{p})$ as $x\downarrow0$, uniformly
for $y$ between $x_1$ and $x_0$.
\\
(2) In
case $T_0f=X$ the product is $1+\oh(x^n)$, for any $n$.
\end{corollary}
\Proof (1) Just use the estimate $f''(x) = \Oh(x^{p-1})$. \\
(2) follows from (1).
\qed

\begin{corollary}\label{corollary-one-all}
If the product (\ref{equation-basic-product})
converges for some $x,\xi>0$, then it converges for
any choice of $x, \xi>0$.
\qed
\end{corollary}

\begin{corollary}Suppose $f$ and $g$ are conjugate in
$\Diffeo([0,\infty))$. Then for any $x>0$ and $\xi>0$ the product
(\ref{equation-basic-product})
converges.
\qed
\end{corollary}

\begin{corollary} The convergence or divergence
of the product (\ref{equation-basic-product}) is not affected if
the functions $f$ and $g$ are replaced by conjugates.
\qed
\end{corollary}

Condition (P) is actually a consequence of Condition (T)
in the non-flat cases:

\begin{proposition}
(1)
If $f'(0)\neq 1$ or $g'(0)\neq 1$, then Condition (P)
    is equivalent to $f'(0)=g'(0)$.
\\
(2) If $f$ and $g$ have conjugate non-identity Taylor series, then
    Condition (P) is satisfied.
\qed
\end{proposition}
\begin{proof}
{(1)}
To prove the first assertion, assume that $f'(0)=\alpha <
1$.  Choose  $\alpha'$ with $\alpha<\alpha'<1$.  Then,  for sufficiently-small
$x$ we have $f(x)<\alpha' x$. Consider the product $\prod\limits_{n=
0}^{\infty} \frac{f'(f^{\circ n}(x))}{\alpha}$. The product
converges if and only if $\sum\limits_{n=0}^{\infty} (1-\frac{f'(f^{\circ
n}(x))}{\alpha})$ converges. But the second derivative of $f$ is bounded
near $0$, so the sum is dominated by a constant times
$\sum\limits_{n=0}^{ \infty}
f^{\circ n}(x)\leq \sum\limits_{n=0}^{
\infty}(\alpha')^n<\infty$, and hence is indeed convergent.

Now  consider the similar product for $g$. The product in
Condition (P) is the quotient of the products if $f'(0)=g'(0)$, and
hence converges as well, and (for the same reason)
it does not converge (to a nonzero limit)
if $f'(0)\neq g'(0)$.

{(2)} By replacing $g$ with a conjugate which has the same
Taylor series as $f$ we reduce to the case in which
$f$ and $g$ have coincident Taylor series.
The result then
follows from the next, more general lemma, which we
will also use again later.
\end{proof}

\begin{lemma}\label{1+O(x^{p-1})}
Let $T_0(f)=T_0(g)=X+bX^{p+1}+\cdots$ (mod $X^{2p+1}$),
where $p\in\N$ and $b\not=0$. Then
$$\prod\limits_{n=0}^{\infty}\frac{f'(f^{\circ n}(x))}{g'(g^{\circ n}
(x))}=1+O(x^{p}).$$
\end{lemma}

\Proof Without loss in generality, we take $b<0$, and write $c=-b$.
We use $C$ for a positive constant that may differ at each
occurrence. We may assume that the $x>0$ under consideration are so
small, that $|f(x)-x+cx^{p+1}|\leq Cx^{p+2}$ and $|Cx|<\frac12 c$.
This means that $c x_n^{p+1}-Cx_n^{p+2}<x_n-x_{n+1}<c x_n^{p+1}+C
x_n^{p+2}$,
where $x_n=f^{\circ n}(x)$. So, for $0<\alpha<1$ between $x$ and $\alpha x$ there
are no more than
$$\frac{(1-\alpha) x}{(c
(\alpha x)^{p+1}- C
x^{p+2})}=\frac{(1-\alpha)\alpha^{-(p+1)}}{cx^{p}(1- \alpha^{-(p+1)}
C x)}
$$ and
no fewer than
$$\frac{(1-\alpha) x}{(cx^{p+1} +Cx^{p+2})}= \frac{1-\alpha}{c
x^{p}(1+Cx)}
$$ points from the $f$-orbit of $x$.

Let us start by reformulating the claim: It is enough to prove
that
$$\log\left(\prod\limits_{n= 0}^\infty \frac{f'(f^{\circ
n}(x))}{g'(g^{\circ n}(x))}\right) =\sum\limits_{ n=0}^\infty
\log\left(\frac{f'(f^{\circ n}(x)) }{g'( g^{\circ n}(x))}\right)
=\Oh(x^{p}).
$$
As $f'(0)=g'(0)=1$ and $|\log(t)|\sim |1-t|$ close to $t=1$, it is
enough to prove that $$\sum\limits_{n=0}^\infty |g'(g^{\circ
n}(x))-f'(f^{\circ n}(x))|=O(x^{p}).$$ Since $T_0f=T_0g$, we may
also assume $x$ is so small that $|f'(x)-g'(x)|<C x^{2p}$. We then
observe that, since  $|x_{n}-x_{n+1}|>(c/2) x_n^{p+1}$, we have the
estimate
\begin{equation}\label{5th_formula_p14}
\sum\limits_{k=0}^\infty x_k^{p+1}\leq
\frac2{c}\sum_{k=0}^\infty |x_k-x_{k+1}| =\frac{2x_0}{c}.
\end{equation}
As $|f'(x)-g'(x)|<C
x^{2p}$ for all $x$ in question, we have
$$\sum\limits_{k=0}^\infty
|f'(x_k)-g'(x_k)|\leq C\sum\limits_{k=0 }^\infty x_k^{2p}
 \leq Cx_0^{p-1}\sum\limits_{k=0 }^\infty x_k^{p+1}
= O(x_0^p),$$
and the estimate can be reduced to estimating $\sum
|g'(x_n)-g'(g^{\circ n}(x))|$. Since $g''=O(x^{p-1})$ we have
$|g'(r)-g'(s)|\leq O(s^{p-1}) |r-s|$ for $r<s$, and it remains to
show that $\sum |f^{\circ n}(x)-g^{\circ n}(x)|< C x$.

Let us now consider only points so close to the origin that
$|f(x)-g(x)|<Cx^{2p+1}$. For those points we have the estimate
$$
\begin{array}{rcl}
|f^{\circ n}(x)-g^{\circ n}(x)|
&\leq& |f^{\circ n}(x)-f(g^{\circ(n-1)}(x))|
+ |f (g^{\circ (n-1)})(x)-g^{\circ n}(x)| \\
&\le&  M_x|f^{\circ(n-1)}(x)-g^{\circ(n-1)}(x)|+ C(g^{\circ (n-1)}(x))^{2p+1}\\
&\le&C(M_x^n+\cdots+1)x^{2p+1},
\end{array}
$$
where $M_x$ is the maximum of $f'$ on the interval $[0,x]$, and thus
can (for small $x$) be estimated from above by $1$ (since $b<0$).
This gives us $|f^{\circ n}(x)-g^{\circ n}(x)|\leq Cn x^{2p+1}$.

Let us consider the first point in the orbit of $x$ with respect to
$f$ which is less than $\alpha x$. Let it be $f^{\circ n_1}(x)$.
Then by the observation at the beginning of the proof, for
$\alpha>\frac12$,
and $x$ sufficiently small,
$n_1<(1-\alpha)C/{x^{p}}$, where the constant
depends only on the Taylor expansion. By the previous paragraph, for
any $k\leq n_1$
$$
|f^{\circ k}(x)-g^{\circ k}(x)| \leq\frac{(1-\alpha)Cx^{2p+1}}{x^p} =
(1-\alpha)C x^{p+1}.
$$

As $|f(y)-y|> \frac{c}2 y^{p+1}$, we see that for a choice of
$\alpha<1$ close enough to $1$, $\frac12(f^{\circ (k+1)}(x)+f^{\circ
k}(x))<g^{\circ k}(x)<\frac12(f^{\circ k}(x)+f^{ \circ (k-1)}(x))$.
Thus the intervals $[f^{\circ k}(x), g^{\circ k}(x)]$ are disjoint,
and $\sum\limits_{k=0}^{n_1} |f^{\circ k}(x)-g^{\circ k}(x)|\leq
(1-\alpha)x+Cx^{p+1}$. On the other hand, in the particular case,
$k=n_1$, if $g^{\circ (n_1+1)}(x)=x^{(1)}$, we have
$\sum\limits_{m=0}^{\infty}|f^{\circ (n_1+1+m)}(x)-f^{\circ
m}(x^{(1)})|\leq \alpha x$, as the sum of lengths of disjoint
intervals. This means that

$$\sum\limits_{n=0}^{\infty} |f^{\circ
n}(x)-g^{\circ n}(x)| $$ $$\le  \sum\limits_{n=0}^{n_1} |f^{\circ
n}(x)-g^{\circ n}(x)|+ \sum\limits_{m=0}^{\infty}|f^{\circ
m}(x^{(1)})-g^{\circ m}(x^{(1)})| +$$
$$\sum\limits_{m=0}^{\infty}|f^{\circ (n_1+1+m)}(x)-f^{\circ
m}(x^{(1)})|+C x^{p+1} $$ $$\leq (1-\alpha)x+
\sum\limits_{m=0}^{\infty}|f^{\circ m}(x^{(1)})-g^{\circ
m}(x^{(1)})|+ \alpha x +Cx^{p+1}\leq $$ $$2x+
\sum\limits_{m=0}^{\infty}|f^{\circ m}(x^{(1)})-g^{\circ
m}(x^{(1)})|.$$

Using this argument inductively we deduce that
$$\sum\limits_{n=0}^{\infty} |f^{\circ
n}(x)-g^{\circ n}(x)| \le 2x+ 2x^{(1)}+\ldots\leq 2
\sum\limits_{j=0}^{\infty} \alpha^j x= Cx,$$ and we are done. \qed

\begin{example}
\end{example}
Notice, that for the particular case $p=1$ this lemma
says that the Condition (P) is satisfied for $f(x)=x+x^2$ and
$g(x)=x+x^2+x^3$. On the other hand, the Taylor series $X+X^2$
and $X+X^2+X^3$ are not conjugate, which shows that the
Condition (P) is strictly weaker than Condition (T) in the non-flat case .

\medskip
We shall see shortly that condition (P) guarantees
the existence of a $C^1$ diffeomorphism conjugating $f$ to $g$.
Thus the existence of a $C^1$ conjugacy is strictly weaker than the
existence of a $C^\infty$ conjugacy.

\medskip
We mention here the observations of Young [Y].  He considered
$C^2$ diffeomorphisms $f$ on $[0,+\infty)$ with
$T_0f=x+ax^2$ (mod$x^3$), and with $a\not=0$.
A result of Szekeres (cf. \cite[Theorem 8.4.5]{KCG}) implies that all such $C^2$ diffeomorphisms
(having $a$ of one sign) are $C^1$-conjugate.
Young showed that they are in fact $C^2$-conjugate.

\section{$[0,+\infty)$: Sufficient Conditions}\label{section-sufficient}

\subsection{The Differential Equation}
Suppose $f,g\in\Diffeo([0,+\infty))$ fix only $0$,
both belong to $S_-$ and satisfy
condition (P).

\medskip
We define
$$ F_{1a}(x) = H_1(f,f;x,a) =
\prod_{n=0}^\infty
\frac{f^\prime(f^{\circ n}(x))}%
{f^\prime(f^{\circ n}(a))}$$
whenever $a,x>0$.  Note that
$$ F_{1a}(x) = \lim_{n\uparrow\infty}
\frac{(f^{\circ n})^\prime(x)}{(f^{\circ n})^\prime(a)}.$$
We define
$$ G_{1\alpha}(\xi) = H_1(g,g;\xi,\alpha) =
\prod_{n=0}^\infty
\frac{g^\prime(g^{\circ n}(\xi))}%
{g^\prime(g^{\circ n}(\alpha))}$$
whenever $\alpha,\xi>0$.
\begin{lemma}\label{lemma-smooth-products}
Fix $a>0$, $\alpha>0$.
The functions $x\mapsto F_{1a}(x)$
and $\xi\mapsto G_{1\alpha}(\xi)$
are
infinitely-differentiable and positive on $(0+\infty)$,
and hence
$$ (x,\xi)\mapsto H_1(x,y) = H_1(a,\alpha)F_{1a}(x)/G_{1\alpha}(\xi)$$
is
infinitely-differentiable and positive on $(0,+\infty)\times (0,+\infty)$.
\end{lemma}
\Proof

It suffices to show that $x\mapsto F_{1a}(x)$
is infinitely-differentiable on $(0,+\infty)$ for each $a>0$.
The argument for $\xi\mapsto G_{1\alpha}(\xi)$ is precisely
analogous.

Fix $a\in (0,+\infty)$. Let $J_a$ denote the closed interval from $0$ to $a$.
Let $a_n= f^{\circ n}(a)$, for all $n\in \Z$.
Let $I_a$ denote the closed interval from $a_1$ to $a$. Let

$$ D_j=\max\limits_{J_a}|f^{(j)}|, \forall j\in \Z.$$
(Note that $\min\limits_{J_a}|f'|=(D_{-1})^{-1}$.)

For $x\in (0,+\infty)$, let $x_n=f^{\circ n}(x)$, for all
$n\in \Z$. For ease of notation, we abbreviate $\frac{d}{dx}f^{\circ n}(x)
=f'(x)f'(x_1)\ldots f'(x_{n-1})$ to $x'_n$, and similarly denote $\frac{d^k}{
dx^k}f^{\circ n}(x)$ by $x^{(k)}_n$. We use $x''_n$ for
$x_n^{(2)}$,
etc.

Before continuing the proof, we pause to note a couple of lemmas that
follow from Lemma \ref{lemma-product-existence}.

In what follows, unless otherwise specified, we use $K$ to denote a
constant that depends at most on $f$, and $a$, and that may be
different at each occurence.

\begin{lemma}
$$K^{-1} |(f^{\circ n})'(a)|\leq |x'_n|\leq K|(f^{\circ n})'(a)|$$
whenever $x\in I_a$.
\end{lemma}
\Proof
$$\frac{(f^{\circ n})'(a)}{(f^{\circ n})'(x)}=\prod\limits_{j=0}^{n-1}
\frac{f'(a_j)}{f'(x_j)},$$ so the result follows from the uniform convergence
of $\prod\limits_{j=0}^{\infty} \frac{f'(a_j)}{f'(x_j)},$ for $x\in I_a$.
\qed

\begin{lemma}
$$|x_n'|\leq K\Big|\frac{x_{n+1}-x_n}{x_1-x_0}\Big|$$
whenever $x\in I_a$.
\end{lemma}
\Proof By the Law of the Mean,
$$\frac{x_{n+1}-x_n}{x_1-x_0}=(f^{\circ n})'(y)$$
for some $y$ between $x$ and $x_1$, so the result follows from a
few applications of the previous lemma.
\qed

\begin{lemma} $|x_1-x|\geq K|a_1-a|$, for all $x\in I_a$.
\end{lemma}
\Proof For $x\in I_a$, $f(a)\le x\le a$,
so $f(x)\le f(a)\le x$, so
$|f(x)-x|=|f(x)-f(a)|+|x-f(a)|\geq (D_{-1})^{-1}|x-a|+|x-a_1|\geq
\min\{1, (D_{-1})^{-1}\}|a-a_1|.$\qed

\noindent
{\bf Proof of Lemma \ref{lemma-smooth-products}.}
It suffices to show that the logarithm
$$\log F_{1a}(x,\xi)=\sum\limits_{n=0}^{+\infty} \{\log f'(x_n)
-\log f'(a_n)\}$$
is infinitely-differentiable.

The term by term derivative with respect to $x$ is the series
$$\sum\limits_{n=0}^{+\infty}\frac{f''(x_n)x_n'}{f'(x_n)},$$
and it will be convenient to denote the $n$-th term by
$$T_n(x)=\frac{f''(x_n)x_n'}{f'(x_n)},$$
and the $n$-th partial sum by
$$S_n(x)=\sum\limits_{j=0}^{n-1}T_j(x).$$

It will suffice to show that for each nonnegative
integer $k$, $S^{(k)}_n(x)$ converges uniformly on $I_a$.

For any smooth function $\rho:(0,+\infty)\to(0,+\infty)$, and $k\in\N$,
let us define $ A_k(\rho)$ as the function
$$ A_k(\rho) = \frac{d^k}{dx^k}\left(
\frac{f''(\rho)\rho'}{f'(\rho)}
\right) -
\frac{f''(\rho)\rho^{(k+1)}}{f'(\rho)}.
$$
Then a straightforward induction
establishes that $A_k(\rho)(x)$ is the sum of $M_k$ terms (where the integer
$M_k$ depends on $k$, but not on $\rho$), each of which is a
finite product
$$\gamma\prod\limits_i f^{(r_j)}(\rho(x))\prod\limits_j
\rho^{(t_j)}(x)\over(f'(\rho(x)))^{k+1},$$
where the coefficients $\gamma$ are fixed integers independent of $f$,
where each $r_i\leq k+2$, each $t_j\leq k$, and at least one $t_j$ is present.

The term $T^{(k)}_n$ takes the form
$$A_k(x_n)+\frac{f''(x_n)x_n^{(k+1)}}{f'(x_n)}$$.

To begin with, we observe that by the last two lemmas
$$\left|\frac{f''(x_n)x'_n}{f'(x_n)}\right|\leq K D_2
D_{-1}|x_{n+1}-x_n| , \qquad\forall x\in I_a$$ hence $\{S_n(x)\}$
itself converges uniformly on $I_a$, with the error in $S_n(x)$
bounded by $K D_2 D_{-1} a_n,$ where $a_n=f^{\circ n}(a)$.

Now we will proceed by induction on $k$, and we
first consider the first derivatives $T_n'(x)$ and note that
$$T_n'(x)= A_1(x_n)+B_1(x_n),$$
where
$$ A_1(x_n)=
\left\{\frac{f'''(x_n)(x'_n)^2}{f'(x_n)} -
 \frac{(f''(x_n)x'_n)^2}{(f'(x_n))^2}\right\}
$$
and
$$ B_1(x_n)=
 \frac{f''(x_n)x''_n}{f'(x_n)}.$$

Estimating each of its terms by its maximum, we see
that $A_1(x_n)$ is dominated by

$$ K^2 (D_3
D_{-1}+D_2^2 D_{-1}^2) |x_{n+1}-x_n|\le K|x_{n+1}-x_n|,
$$ for a (different) constant $K$.

A calculation yields $x''_n=x'_n S_n$, so the term $B_1(x_n)$ is
dominated  by
$$K^2D_2D_{-1}|x_{n+1}-x_n|,$$
and we conclude that $S'_n(x)$ also converges uniformly on $I_a$, with
error bounded by $Kx_n$.

We also observe that $|S'_n(x)|\leq Kx$.

Now we formulate an {\bf induction hypothesis} $P_k$:

{\it There exist a constant $K$, depending only on $f$, $a$, and
$k$, such that

$(a)$ for $0\leq j\leq k-1$, and each $n\ge0$,
$$|T^{(j)}_n(x)|\leq K|x_{n+1}-x_n|
\hbox{, and } |S^{(j)}_n(x)|\leq K\hbox{, and}$$

$(b)$ for $1\leq j\leq k$,
$$|x^{(j)}_n|\leq K|x_{n+1}-x_n|.$$}

We have established $P_2$.

Suppose $P_k$ holds, for some $k\geq 2$. Differentiating the formula
$x''_n=x'_n S_n$ $k-1$ times we get
$$ x^{(k+1)}_n=\sum\limits_{j=0}^{k-1}{\binom{k-1}{j}}x_n^{(j+1)}S_n^{
(k-j-1)}$$
so conditions (a) and (b) of the hypothesis yield
$$|x^{(k+1)}_n|\leq \sum\limits_{j=0}^{k-1}{\binom{k-1}{j}}K|x_{n+1} - x_n|=
K|x_{n+1}-x_n|$$
(with a new $K$), and condition (b) of $P_{k+1}$ is proven.

Condition (a) then follows because of the form of $T^{(k)}_n$.

Thus, by induction, $P_k$ holds for each $k\geq 2$.

Thus $S^{(k)}_n=\sum T_m^{(k)}$ converges uniformly for all $k$, and
Lemma \ref{lemma-smooth-products} is proved.

\qed

\medskip
We now consider the three-parameter initial-value
problem
$D_1(a,\alpha,\lambda)$
(cf. equation (\ref{equation-IVP-D1})).

It follows from Lemma \ref{lemma-smooth-products}
and standard results about
ordinary differential equations \cite[p. 22]{BDP}
 that problem
$D_1(a,\alpha,\lambda)$ has a unique infinitely-differen\-tiable
solution $\phi(a,\alpha,\lambda;x)$ near $x=a$ whenever
(P) holds,
$a,\alpha>0$, and $\lambda>0$. Obviously,
the solution is a strictly-increasing function of $x$
and its domain is an open subinterval of $(0,+\infty)$, containing $a$.

Note that
$$ \prod_{j=n}^mf^\prime(x_j) = (f^{\circ(m-n)})^\prime(x_n)
\approx{x_{m+1}-x_m\over x_{n+1}-x_n},$$ so the product tends to $0$
as $m\to+\infty$. It follows that the product $H_1(x,y)$ does not
extend continuously to the closed quadrant $[0,+\infty)\times
[0,+\infty)$, nor even to the corner $(0,0)$, so there is no point
in considering the differential equations at the endpoint. In fact,
a moment's thought reveals that $H_1(x,y)$ tends to $\infty$ as
$x\to 0$ for fixed $y>0$, and tends to $0$ as $y\to 0$ for fixed
$x>0$, so all positive numbers may be obtained as limits of
$H_1(x,y)$ for suitable approach to $(0,0)$ from inside $J\times J$.

\begin{lemma} Assume $f,g\in S_-$ and condition (P) holds.
Then for each $a,\alpha>0$ and each $\lambda>0$
the domain of the solution to problem $D_1(a,\alpha,\lambda)$ is $(0,+\infty)$.
\end{lemma}
\begin{proof}
The domain $U$ of the solution $\phi$
is a nonempty connected subset of $J=(0,+\infty)$,
and one sees easily that it is open. In fact, if
either end (say $c$) of $U$ lies
inside $J$,
then since $\phi'$ is bounded
on $(a,c)$, $\phi$ has a continuous extension to
$c$, and by uniqueness
$\phi$ is the solution to
$D_1(c,\phi(c),\lambda)$, so
$\phi$ extends to a neighbourhood of $c$, a contradiction. Thus
$U=J$, and $\phi$ conjugates $f$ to $g$ on the whole of $J$.
\end{proof}

The following lemma reformulates the information in the proof of
Lemma \ref{lemma-basic-product}
in other language:
\begin{lemma}
Suppose $f,g,h\in \Diffeo([0,+\infty))$ and $f=g^h$.
Then
$\phi=h$ is the solution
to problem $D_1(a,h(a),h^\prime(0))$, whenever
$a>0$.
\qed
\end{lemma}

Not all solutions to the
initial-value problems $D_1(a,\alpha,\lambda)$ will be conjugating maps.
For a start, we would need to ensure the condition
$\phi(f(a))=g(\alpha)$.
This leads us to the following:

\begin{lemma}\label{lemma-technical}
Assume condition (P), with $f,g\in S_-$. Then
for each $a>0$, and each $\alpha>0$, there exists $\lambda>0$ such that
the solution $\phi(a,\alpha,\lambda)$
to problem $D_1(a,\alpha,\lambda)$ has
$\phi(f(a))=g(\alpha)$.
\end{lemma}
\Proof Given $a$ and $\alpha$, we could start by trying $\lambda=1$.
If the solution $\phi_1$ to $D_1(a,\alpha,1)$ has
$\phi_1(f(a))=g(\alpha)$, we take $\lambda=1$ and are done. If
$\phi_1(f(a))<g(\alpha)$, then decreasing $\lambda$ eventually
reduces $\phi^\prime$ to very small values on the interval
$[a,f(a)]$, and hence pulls $\phi(f(a))$ above $g(\alpha)$. Thus,
since $\phi(f(a))$ varies continuously with $\lambda$, there exists
some $\lambda$ with $\phi(a,\alpha,\lambda)(f(a))=g(\alpha)$, and we
are done. If $\phi_1(f(a))>g(\alpha)$, then we can attain a similar
result by increasing $\lambda$ instead, because this increases
$\phi^\prime$ to very large values. More precisely,
$H_1(x,y)$ is bounded below by a positive constant, say $\kappa$,
 on $[f(a),a]\times[g(\alpha),\alpha]$, so if $\phi(f(a))\ge g(\alpha)$,
then $\phi'>\kappa\lambda$ on $[f(a),a]$, hence
$$ \alpha-g(\alpha) \ge \phi(a)-\phi(f(a)) \ge \lambda\kappa(a-f(a)),$$
which is impossible for large $\lambda$. Thus for large enough
$\lambda$, we have $\phi(f(a))<g(\alpha)$, so another application
of the intermediate value theorem tells us that there
exists some $\lambda$ with $\phi(a,\alpha,\lambda)(f(a))=g(\alpha)$.
\qed

Now we proceed to show that the solution $\phi$
of Lemma (\ref{lemma-technical})
conjugates $f$ to $g$ on
$(0,+\infty)$.

\begin{lemma}\label{lemma-inequality}
Suppose that $u$ is a differentiable real-valued function
on an open interval $U$,
and for some
constant $\kappa>0$ we have
$$ |u^\prime(x)| \le \kappa\cdot|u(x)|,\quad\forall x\in U. $$
Suppose that $u$ has a zero in $U$. Then
$u$ is identically zero on $U$.
\end{lemma}
\Proof The set $Z=u^{\circ-1}(0)$ of zeros
of $u$ in $U$ is relatively-closed, and nonempty,
so it suffices to show that it is open.
Fix $a\in Z$, and choose $\epsilon>0$
so that $a\pm\epsilon\in U$
and $\epsilon\kappa<1$.
Let $M$ be the maximum of $|u|$
on the closed interval
$J=[a-\epsilon,a+\epsilon]$.

If $M>0$, then choose $b\in J$
with $|u(b)|=M$. By the Law of the Mean,
we may choose $c$ between $a$ and $b$
with $|u(b)|=
|u^\prime(c)|\cdot|b-a|$.
But then
$$ M=|u(b)|\le \kappa M\cdot\epsilon < M,$$
which is impossible.

Thus  $M=0$, so $a$ is an interior point of $Z$.

Thus $Z$ is open, and we are done.\qed

\begin{lemma}\label{lemma-conjugation}
Suppose (P). Fix $a,\alpha>0$. Choose $\lambda>0$ such that
the solution $h(x)=\phi(a,\alpha,\lambda;x)$ to problem
$D_1(a,\alpha,\lambda)$ has $h(f(a))=g(\alpha)$. Then the domain of
the solution is $J=(0,+\infty)$, $h$ maps $J$ onto $J$ and $g\circ
h=h\circ f$ on $J$.
\end{lemma}

\Proof
We establish that on each compact subinterval
of $J$ we have an inequality $|u^\prime|\le\kappa\cdot|u|$,
where
$$ u(x) = g(\phi(x))-\phi(f(x)).$$

In detail, one calculates (by fiddling with products) that
$$ u^\prime(x) =
\lambda \{H_1(x,\phi(x))\cdot
g^\prime(\phi(x))-H_1(x,g^{-1}(\phi(f(x)))) \cdot
g^\prime(g^{-1}(\phi(f(x))))\},$$ and (using the Law of the Mean)
estimates this (on a compact subinterval of $J$) by
$$ \kappa_1\times |g^{-1}(\phi(f(x)))-\phi(x)| $$
$$ \le \kappa_2|\phi(f(x))-g(\phi(x))| = \kappa_2|u(x)|.$$

Then we apply Lemma \ref{lemma-inequality} and the fact that
$u(a)=g(\alpha)-\phi(f(a))=0$. This tells us that $u(x)=0$ on the
domain of $\phi$, which is $J$.
\qed

These results tell us that the initial-value problem
together with the conjugation equation at one
point are enough to guarantee the conjugation equation
on the whole interval $J=(0,+\infty)$.

\begin{lemma}\label{lemma-one-sided}
Suppose condition (P) holds.
If $\phi:[0,+\infty)\to [0,+\infty)$ satisfies
$$ \left.
\begin{array}{rcl}
\phi(f(x))&=&g(\phi(x)),\\
\phi^\prime(x) &=& H_1(x,\phi(x))\lambda
\end{array}\right\}\qquad\forall x\in J,
$$
then $\lim_{x\to d}\phi^{\prime}(x)=\lambda$ and $\phi$ has a
one-sided derivative at $d$, equal to $\lambda$.
\end{lemma}
\Proof
Fix some $a\in J$ and denote $I_a=[f(a),a]$.

For fixed $x\in I_a$, letting $x_n=f^{\circ n}(x)$, we have
$$\begin{array}{rcl}
\phi(f^{\circ n}(x)) &=& g^{\circ n}(\phi(x)),\\
\phi^\prime(x_n)\cdot x_n^\prime &=&
(g^{\circ n})^{\prime}(\phi(x))
\cdot\phi^\prime(x),\\
\phi^\prime(x_n) &=&
\prod_{j=0}^{n-1}\left(
{\displaystyle g^\prime(g^{\circ j}(\phi(x))\over\displaystyle  f^\prime(x_j)}
\right)\cdot\phi^\prime(x),
\end{array}.
$$
Since the
product converges to $H_1(x,\phi(x))^{-1}$,
the right-hand side converges to $\lambda$, so the
derivative $\phi^\prime$ extends continuously from $J$ to $0$ if
$\phi$ is given the value $0$ there. This is enough to force the
rest of the conclusions. \qed

Finally, we show that the $\lambda$ is unique:

\begin{lemma}\label{lemma-exactly-one}
Suppose conditions (P) holds.
Then, for each given $a,\alpha\in J$,
there is {\bf exactly one} $\lambda>0$ for which the solution $\phi=h$ to problem
$D_1(a,\alpha,\lambda)$ has $h(f(a))=g(\alpha)$.
\end{lemma}
\Proof
Suppose this fails, and there are $\lambda_1<\lambda_2$
such that the solutions $\phi_i$ to problems
$D_1(a,\alpha,\lambda_i)$ ($i=1,2$)
both have $\phi_i(f(a))=g(\phi_i(a))$.

Then by Lemma \ref{lemma-conjugation} both solutions have
$\phi_i(f(x))=g(\phi_i(x))$ on $J$, both map
$J$ onto $J$, and both derivatives extend
continuously to $0$.

Since, initially,
$\phi_1(a)=\phi_2(a)$ and
$\phi_1^\prime(a)<\phi_2^\prime(a)$,
we have $\phi_1(x)>\phi_2(x)$ for some
distance to the left of $a$.
Since $\phi_1(0)=\phi_2(0) (=0)$, there
exists a first point $e<a$ at which $\phi_1(e)=\phi_2(e)$.
Just to the right of $e$, we have
$\phi_1(x)>\phi_2(x)$, and hence
$\phi_1^\prime(e)\ge\phi_2^\prime(e)$.
But this contradicts the differential equation, because
(since $\phi_1(e)=\phi_2(e)$) we have
$$ \phi_1^\prime(e) = \lambda_1H(e,\phi_1(e))
<\lambda_2H(e,\phi_2(e))=\phi_2^\prime(e).$$
This contradiction establishes the result.
\qed

\medskip
At this stage, we have completed the proof of Theorem \ref{theorem-main}.

\begin{corollary} Suppose Conditions (P) holds.  Then
there is precisely a one-parameter
family of $C^1$ conjugations from $f$ to $g$
on $[0,+\infty)$.
\end{corollary}
\Proof In fact, if we fix $a$, there is precisely one
conjugation $\phi=\Phi_+(a,\alpha)$
for each $\alpha\in(0,+\infty)$.
\qed

Thus there is {\em at most} a one-parameter family of $C^\infty$
conjugations from $f$ to $g$.
(One could recover Kopell's Lemma (cf. Subsection \ref{subsection-remarks-C})
from this. However, it can be proved directly without
all this apparatus (cf. \cite[4.1.1]{N2}). One should also
remark that the Corollary may be obtained directly from
Kopell's Lemma, and holds for $C^2$ conjugations.)

\subsection{Remarks about $\phi^\prime(0)$}
Assume Conditions (P) holds.

If $0$ is a hyperbolic point for $f$, then the family of conjugating maps is parametrised by the multiplier
at $0$. This is so, because $f$ is conjugate to $\lambda\cdot$
for $\lambda=f'(0)$, which has centraliser consisting
of all the maps $\mu\cdot$, with $\mu\not=0$. If two maps
$h$ and $k$ that conjugate $\lambda\cdot$ to $g$ have the same
multiplier, then $k^{-1}h$ commutes with $\lambda\cdot$
and has multiplier $1$, and hence $h'(0)=k'(0)$.

If $f'(0)=1$, but $f-x$ is not flat at $0$, then
it follows from a theorem of Lubin's theorem about the centraliser of a formal
power series \cite{L, OF2}
that the centraliser of $f$
in $\Diffeo([0,+\infty))$ consists of maps that have derivative 1
at $0$, and for general $g$ the family of conjugating maps
from $f$ to $g$ is a coset of this centraliser. Thus
all the diffeomorphic conjugations of $f$ to $g$
have the same derivative at $0$.

This does not tell us anything about merely $C^1$ conjugations,
nor about what
happens when $f-x$ is flat at $0$, but it is possible to see that
again the conjugating $C^1$ maps all have the same derivative
at $0$.  The essential point is the following, which can be proved
more simply now than Lemma \ref{1+O(x^{p-1})}:
\begin{proposition}\label{proposition-4.12}
Suppose $f\in S_-$, $f'(0)=1$, and $\phi$ is a $C^1$ diffeomorphism
of $[0,+\infty)$, commuting with $f$.  Then $\phi'(0)=1$.
\end{proposition}
\Proof  Fix $a>0$, and let $\alpha=\phi(a)$.  Then
$\phi$ is $\Phi_+(a, \alpha)$.  Let $a_k=f^{\circ k}(a)$
whenever $k\in\Z$. There is a unique $k$ such that
$$ a_{k+1} \le \alpha < a_{k}.$$
So at $a$, $\phi$ lies between $f^{\circ k}$
and $f^{\circ(k+1)}$.

If $\phi(a)=f^{\circ k}(a)$, then
by
Lemma
\ref{lemma-exactly-one},
$\phi$
coincides with $f^{\circ k}$ on $J$,
and hence has derivative $1$ at $0$,
and we are done.

Otherwise,
\ref{lemma-exactly-one} tells us that
$\phi$ never has the same value as
$f^{\circ k}$ or
$f^{\circ (k+1)}$  at any point, so
its graph lies sandwiched between
their graphs.

Thus, since $f(0)=\phi(0)=0$,
$$
f^{\circ k}(x)
> \phi(x) - \phi(0)>
f^{\circ(k+1)}(x)-f^{\circ (k+1)}(0)$$
for all $x>0$, and hence, dividing by $x$
and taking limits we get $\phi'(0)=1$.\qed

We remark that this result becomes trivial for $C^2$ conjugations.
If we assume that the conjugating map is
$C^\infty$ to $0$, then Corollary \ref{corollary-flat-ff-product}
provides a much easier way to a stronger conclusion:

\begin{proposition}{\bf (Kopell)}\label{proposition-flat-centraliser}
If $\phi\in\Diffeo([0,+\infty))$ commutes with $f$,
and $f$ is flat at $0$, then so is $\phi$.
\end{proposition}
\Proof From Corollary \ref{corollary-flat-ff-product},
$\phi(x)-x$ tends to
zero more rapidly than any power of $x$,
and hence {\em given that $\phi(x)-x$ is smooth},
all its derivatives vanish at $0$. \qed

\subsection{Remark about Centralisers}\label{subsection-remarks-C}
The special case $f=g$ of the foregoing
corresponds to results of Kopell \cite[pp. 167-71]{K}
about centralisers. Indeed, Kopell made use of
the $f=g$ version of the differential equation
of problem $D_1$ in order to obtain her results.
See also \cite[Section 8.6, pp. 353-5]{KCG}.
(We have not seen the differential equation
for general $f$ and $g$ used in the literature.)

The elements of the centraliser $C_f$ of
$f$ in $\Diffeo([0,+\infty))$ (where
$f$ fixes only $0$)
are exactly the elements that conjugate $f$
to $f$, so applying the foregoing to the case
$g=f$, we have Kopell's result that the centraliser
is at most a one-parameter group.
The centraliser is never trivial, since
it has all iterates $f^{\circ n}$ ($n\in\Z$)
as elements.  However, it may fail to be connected.
Sergeraert \cite{SE}
gave an example in which
$f$ has no smooth compositional square root,
and hence its centraliser is discrete.

Sergeraert also gave a useful sufficient condition
for the centraliser of an element $f\in S_-$
to be connected.  His condition
is the existence of constants $\kappa>0$ and $\delta>0$ such that
$$  \sup_{0\le y\le x}(y-f(y)) \le \kappa (x-f(x)),$$
whenever $0<x<\delta$.  In particular, it always
works if $x-f(x)$ is monotone.

The homomorphism $h\mapsto h^\prime(0)$
maps the centraliser of a given $f$ to
a multiplicative subgroup of $(0,+\infty)$,
but (as we've seen) the subgroup in question is just $\{1\}$,
as soon as $f'(0)=1$.

In a rather similar way, the homomorphism
$\Pi:h\mapsto T_0h$ maps $C_f$ onto
a subgroup of the group
of invertible formal power series,
and the image must have $T_0f$ as an element.

We have seen in Proposition \ref{proposition-flat-centraliser}
that
if $f-x$ is flat at $0$,
then all elements
of its centraliser have the same property, so $\Pi$ is trivial.

Generally, the image of $C_f$ under $\Pi$ is a subgroup of the
centraliser of $T_0f$ in the power series group.
In case $T_0f=X$ mod$X^{p+1}$ but
$T_0f\not=X$ mod$X^{p+2}$, it is a purely algebraic
fact (cf. \cite[p. 170]{K}, \cite{L} or \cite[p. 355ff]{KCG}) that the latter
centraliser is a one-parameter group,
and indeed the map to the coefficient of $X^{p+1}$
is an isomorphism to $(\R,+)$.

\medskip
It is interesting to note in passing that
the differential equation provides a way to construct smooth compositional
$k$-th roots of a diffeomorphism $f\in S^-$ of $[0,+\infty)$
that has a connected centraliser:
One takes $f=g$, fixes $a>0$,
and considers the initial-value problem
 $D_1(a, \alpha, \Lambda_+(a,\alpha))$
for $\alpha$ between $a$ and $f(a)$. The solution
$\phi_\alpha$ that has $\phi^{\circ k}(a)=f(a)$
is the desired root.
Since $\phi^{\circ k}(a)$ moves continuously and
monotonically away from
$a$ as $\alpha$ moves towards $f(a)$
from $a$, and passes $f(a)$ before $\alpha$ reaches
$f(a)$, there must exist a unique $\alpha$
with the above property.

In general, if $f$ does not have roots in the diffeomorphism group,
this will construct $C^1$ roots.

\subsection{Sufficiency of (P) and (T): Counterexample}
The conditions (P) and (T) together are not
sufficient for $C^\infty$ conjugacy, and
the following example will demonstrate this.

We have noted that in the non-flat case
the existence of a $C^1$-conjugacy is strictly weaker than
the existence of a $C^\infty$ conjugacy.
The example will also show that it is also weaker in the flat case.

\begin{example}
\end{example}
Consider
the diffeomorphisms of $[0,+\infty)$ defined
on the interior by
$$\begin{array}{rcl}
f(x) &=& x + e^{-1/x^2},\\
\phi(x) &=& x + x^{3/2},\\
g &=& f^\phi.
\end{array}
$$
One finds that $\phi$ is only $C^1$ on $[0,+\infty)$, but
that $f$ and $g$ are smooth:
In fact, letting $\psi=\phi^{\circ -1}$, we calculate

$$\psi'(\phi)\phi'=1,$$
\begin{equation}\label{equation-2deriv}
\psi'(\phi)\phi''+\psi''(\phi)(\phi')^2=0.
\end{equation}
Thus
$$g'=\psi'(f\circ \phi)f'(\phi)\phi',$$
$$g''=\psi'(f\circ \phi)f'(\phi)\phi''+\psi'(f\circ
\phi)f''(\phi)(\phi')^2+\psi''(f\circ \phi)\{f'(\phi)\phi'\}^2.$$

The second term in the expression for $g''$ is continuous, and the
other two add to
\begin{equation}\label{equation-2'deriv} f'(\phi)\{\psi'(f\circ
\phi)\phi''+\psi''(f\circ \phi)f'(\phi)(\phi')^2\}.
\end{equation}
The only problem is to see continuity at $0$, and the point is that
for small positive $x$ we have $\phi'(x)\approx 1$, $\psi'(x)\approx
1$,
$$\phi^{(k)}(x)=O(x^{\frac32-k}), \forall k\geq 2$$
and for some sequence of integers $p_k$,
$$\psi^{(k)}(x)=O(x^{-p_k}), \forall k\geq 2$$
(as is verified inductively
by differentiating (\ref{equation-2deriv})).
 Thus, since $f(x)-x$ is flat at $0$,
$f'(\phi(x))$ may be replaced by $1$ and $f\circ \phi$ by $\phi$, in
the expression (\ref{equation-2'deriv}), with an error that is
$O(x^N)$ for all $N\in \N$. But when this is done we just get $0$,
by (\ref{equation-2deriv}), so $g''\rightarrow 0$ as $x\rightarrow
0$.

It now becomes clear that when we continue to differentiate $g$, and
express $g^{(k)}$ in terms of $\psi, f$, and $\phi$, we get an
expression involving derivatives of $\psi$ (at $f\circ \phi$), $f$
(at $\phi$), and $\phi$, and that when $f$ is replaced by $\ONE$ in
this expression we get zero (the $k$-th derivative of $\psi\circ
\phi$). Moreover, for small $x$, the error involved in replacing
$f(\phi)$ by $\phi$, $f'(\phi)$ by $1$, $f''(\phi)$ by $0$, and all
higher derivatives $f^{(k)}(\phi)$ by $0$, is $O(x^N)$ for all $N$.
Thus $g^{(k)}\rightarrow 0$ as $x\rightarrow 0$ for all $k\geq 3$,
as well. It follows that $g$ is $C^{\infty}$, and $g(x)-x$ is flat
at $0$, as required.

Now any other $C^1$ conjugation of $f$ and $g$ will differ from
$\phi$ by composition with an element of the centralizer of $f$.
Since $f(x)-x$ is monotone, it satisfies Sergeraert's condition
\cite[p.259, Th.3.1]{SE},
and hence the centralizer of $f$ consists of
$C^{\infty}$ diffeomorphisms, and hence no conjugation of $f$ to $g$
is better than $C^1$.

This shows that Conditions (P) and (T) are not sufficient, by themselves,
to guarantee conjugacy,
in general.

\medskip
\noindent
{\bf Question.}
Since not all $C^1$ conjugacies between  a given $f$ and $g$ belonging to
$\Diffeo([0,\infty))$ are $C^\infty$ to zero, it would be interesting to know
whether or not the set of parameters $\alpha$ for which
the solution
$\Phi_+(a,\alpha)$
is $C^\infty$ to zero is
always a relatively closed subset of $(0,\infty)$. We were not able to resolve
this question\footnote{Added in proof: The question has now been answered by H. Eynard-Bontemps, who constructed \cite{EY} an $f\in\Diffeo([0,+\infty))$ whose $C^2$-centraliser is an uncountable  proper subset of its $C^1$ centraliser. Taking $g=f$, this implies that the set of parameters $\alpha$ for which $\Phi_+(1,\alpha)$ is $\C^\infty$ is not relatively-closed in $(0,+\infty)$.}.

\section{Compact Intervals}\label{section-compact}

Recall from Subsection \ref{subsection-compact}
that in the context of a compact interval
the meaning of
Condition (P) must now be modified, and that it
now involves a two-sided product.

\subsection{The sign condition and Condition (P)}
Let $I=[d,c]$ be a compact interval with interior $J$.
Let $f$ and $g$ belong to $\Diffeo^+(I)$.
The sign condition
is, as before, necessary for conjugacy of $f$ and $g$.

Applying Lemma \ref{lemma-product-existence}
and Corollary \ref{corollary-one-all}
to the inverse maps $f^{\circ-1}$
and $g^{\circ-1}$ on the half-open
interval $J\cup\{c\}$, we see that similar results
hold for the products
$$ \prod_{n=1}^\infty{ f^\prime(x_{-n})\over
f^\prime(y_{-n})}
\hbox{ and }
 \prod_{n=1}^\infty{ f^\prime(f^{\circ-n}(x))\over
g^\prime(g^{\circ-n}(\xi))}.
$$
Thus we obtain:
\begin{lemma}\label{lemma-product-compact1}
Suppose $I$ is compact and
$f$ and $g$ are conjugate in $\Diffeo^+(I)$.
Then:
\begin{equation}\label{equation-product-forward}
 \prod_{n=0}^\infty{ f^\prime(f^{\circ n}(x))\over
g^\prime(g^{\circ n}(\xi))}
\end{equation}
converges for each $x,\xi\in J$, and
\begin{equation}\label{equation-product-backward}
 \prod_{n=1}^\infty{ f^\prime(f^{\circ-n}(x))\over
g^\prime(g^{\circ-n}(\xi))}
\end{equation}
converges for each $x,\xi\in J$.
\end{lemma}
\qed
The proof of Lemma \ref{lemma-basic-product} gives:
\begin{corollary}
If $f=g^h$, then
the two-sided product
$$ \prod_{n=-\infty}^\infty{ f^\prime(f^{\circ n}(x))\over
g^\prime(g^{\circ n}(h(x)))} $$
is independent of $x\in J$, and equals the ratio
$h^\prime(c)/h^\prime(d)$ of the derivatives of
the conjugating map at the ends.
\end{corollary}
\qed

\subsection{The Differential Equation}
Suppose $f,g\in\Diffeo^+(I)$ satisfy the sign condition
and condition (P).

\Definition
We define
$$ F_{1a}(x) = H_1(f,f,x,a)=
\prod_{n=0}^\infty {f^\prime(f^{\circ n}(x))\over
f^\prime(f^{\circ n}(a))}$$
$$ G_{1\alpha}(\xi) = H_1(g,g,\xi,\alpha)=
\prod_{n=0}^\infty
 {g^\prime(g^{\circ n}(\xi)) \over
g^\prime(g^{\circ n}(\alpha))}$$
$$ F_{2a}(x) =H_2(f,f,a,x)=
\prod_{n=1}^\infty{ f^\prime(f^{\circ-n}(x))\over
f^\prime(f^{\circ-n}(a))}$$
$$ G_{2\alpha}(\xi) = H_2(g,g,\alpha,\xi)=
\prod_{n=1}^\infty{ g^\prime(g^{\circ-n}(\xi))\over
g^\prime(g^{\circ-n}(\alpha))}$$
whenever $x,\xi,a,\alpha\in J$.

Applying Lemma \ref{lemma-smooth-products}
to the original maps and to their inverses,
we obtain:

\begin{lemma}\label{lemma-smooth-products-compact}
Let the sign condition and (P) hold, and fix $a,\alpha\in J$.
Then \\
(1) The functions
$F_{1a}$, $G_{1\alpha}$, $F_{2a}$, and $G_{2\alpha}$
are infinitely-differentiable and positive on $J$,
and hence
\\
(2) $H_1 (x,\xi) =
H_1(a,\alpha)\cdot F_{1a}(x)/G_{1\alpha}(\xi)$ and
$H_2(x,\xi) =
H_2(a,\alpha)\cdot G_{2\alpha}(\xi)/F_{2a}(x)$
are infinitely-differentiable and positive on $J\times J$.
\end{lemma}
\qed

We now consider {\em two} three-parameter initial-value
problems (2) and (4):
$$ D_1(a,\alpha,\lambda):\qquad
\left\{
\begin{array}{rcl}
 {\displaystyle d\phi\over \displaystyle dx} &=&
H_1(x, \phi(x))\lambda,\\
\phi(a)&=&\alpha
\end{array}
\right.$$
$$ D_2(a,\alpha,\mu):\qquad
\left\{
\begin{array}{rcl}
{\displaystyle d\phi\over\displaystyle dx} &=&
{H_2(x, \phi(x))\mu},\\
\phi(a)&=&\alpha.
\end{array}
\right.$$

Applying Lemma \ref{lemma-product-existence} and Corollary 5.2
twice, we obtain:

\begin{lemma}
Suppose $f,g,h\in \Diffeo^+(I)$ and $f=g^h$. Then\\
(1) The restriction $\phi=h|J$ is the solution
to problem $D_1(a,h(a),h^\prime(d))$, whenever
$a\in J$; \\
(2) the same
$\phi$ is the solution to problem
$D_2(a,h(a),h^\prime(c)^{-1})$,
whenever $a\in J$;
\\
(3)
the function $a\mapsto H(a,h(a))$ is constant on $J$, equal to
$h^\prime(c)/h^\prime(d)$.
\end{lemma}
\qed

To characterise the existence of a conjugating $h$, we need
to formulate the conditions of this lemma in a way that does
not refer explicitly to $h$.
As before, we can do this by using the differential equations.
The following is a consequence of the series of lemmas
from the last section.

\begin{lemma}
Assume the sign condition and (P).\\
(1) For each $a\in J$,
and each $\alpha\in J$, there exists a unique $\lambda>0$ such that
the unique solution $\phi$ to problem $D_1(a,\alpha,\lambda)$ has
$\phi(f(a))=g(\alpha)$.  This $\phi$ is a bijection of
$J$ onto $J$, and has a one-sided derivative at $d$,
with
$$\lim_{x\to d}\phi'(x)=\phi'(d)=\lambda.$$
\\
(2)
For each  $a\in J$,
and each $\alpha\in J$, there exists a unique $\mu>0$ such that
the unique solution $\psi$ to problem $D_2(a,\alpha,\mu)$ has
$\psi(f(a))=g(\alpha)$.  This $\psi$ is a bijection of $J$ onto $J$,
and has a one-sided derivative at $c$, with
$$ \lim_{x\to c}\psi'(x)=\psi'(c)=\mu.$$
\end{lemma}
\qed

So {\em either} initial-value problem
together with the conjugation equation at one
point are enough to guarantee the conjugation equation
on the whole interior $J$.

With the notation of the last lemma,
recall that we denote the unique
$\lambda$
of part (1) by $\Lambda_+(a,\alpha)$, and the
corresponding $\phi(x)$ by $\Phi_+(a,\alpha;x)$.
Similarly, we
denote the $\mu$ of part (2)  by $\Lambda_-(a,\alpha)$ and
the $\psi$ by $\Phi_-(a,\alpha;x)$.

\subsection{Extending a conjugation to ends of $I$}
Assuming the sign condition and (P),
we consider the following condition:

\medskip
\noindent{\bf Condition (E):}\\
{\sl There exist
$a, \alpha\in J$,
for which the solution
$\phi=\Phi_+(a,\alpha;\cdot)$
has a $C^\infty$ extension
to $I$ (and hence agrees with the
solution $\Phi_-(a,\alpha;\cdot)$). }

\medskip
It is clearly equivalent to
replace  \lq\lq there exist $a,\alpha$" by
\lq\lq for each $a$ there exists $\alpha$" in
the formulation of Condition (E).

We note the following.

\begin{lemma} Suppose the sign condition, (P), and (E).
Then if $h=\phi$ is a solution to problem $D_1(a,\alpha,\lambda)$ with $g(\alpha)=\phi(f(a))$
and extends smoothly to the ends of $J$, it follows that
\\
(1) $h^\prime(d) = \lambda$;
\\
(2) $f=g^h$ on clos$J$;
\\
(3) $h$ is a solution to problem $D_2(a,\alpha,h^\prime(c))$;
\\
(4) $h^\prime(c) = H(a,\alpha) h^\prime(d)$;
\end{lemma}
\qed

The proof of Theorem \ref{theorem-compact}
is now complete. \qed

\subsection{Functional Moduli}\label{subsection-modulus}
Robbin \cite[p.424]{RO} described the solution to the
conjugacy problem on a compact interval, subject to the sign condition,
in the case when both ends are hyperbolic fixed
points, i.e.
when $f^\prime(c)\not=1\not=f^\prime(d)$. (See also
\cite{B} and \cite[Chapter 2]{KH}.)

In that case, condition (P) reduces to
the two equations $f^\prime(c)=g^\prime(c)$ and
$f^\prime(d)=g^\prime(d)$, and
conjugacy may be characterised in terms of a modulus.  Robbin's modulus is
a diffeomorphism of $(0,+\infty)$.  He constructs the
modulus for $f$ by linearizing the restrictions of $f$
to $\{c\}\cup J$ and to $J\cup\{d\}$, i.e. choosing
the (unique) $\alpha_f:\{c\}\cup J\to[0,+\infty)$ and
$\beta_f:J\cup\{d\}\to[0,+\infty)$ such that
$$ \begin{array}{rcl}
\alpha(f(x)) &=& f^\prime(c)\cdot\alpha(x),\\
\beta(f(x)) &=& f^\prime(d)\cdot\beta(x),
\end{array}
$$
whenever $x\in J$, and $\alpha^\prime(c)=\beta^\prime(d)=1$. His
modulus is $\gamma_f=\beta\circ\alpha^{\circ-1}$.  The two maps $f$
and $g$ are conjugate if and only if $\gamma_f=\gamma_g$.  Thus
Robbin's modulus serves to label the elements of the
uncountable family of conjugacy classes
with respect to the group
$\Diffeo^+(I)$ on the compact interval
into which each single
conjugacy class with respect to the group
$\Diffeo^+(J\cup\{c\})$
on one of the half-open intervals splits.

One could try to construct an invariant composed
of Taylor series conjugacy classes and
a modulus, for the general non-flat cases.
In fact,
Young \cite{Y} has shown that the {\em conventional multiplier}
introduced by Afraimo\-vitch, Liu and Young \cite{ALY}
can be used to make a modulus for the
\lq\lq saddle-node" case (in which $f(x)-f(p)$
vanishes to precisely second order at the ends
$p$ of $I$.  We expect that this works
as soon as $f-x$ not flat at either end.

Afraimovitch {\it et al.}, associated to suitable diffeomorphism
$f\in S$ and any fixed $a\in J$ the functions
$$ u_+(x) = \lim_{n\uparrow+\infty}{f^{\circ n}(x)-f^{\circ n}(a)
\over f^{\circ (n+1)}(a)-f^{\circ n}(a)
}$$
and
$$ u_-(x) = \lim_{n\downarrow-\infty}{f^{\circ n}(x)-f^{\circ n}(a)
\over f^{\circ (n+1)}(a)-f^{\circ n}(a)
},$$
defined for each $x\in J$.  They showed that the limits
exist when $f-x$ is not flat at $0$.
Also, in that case,
each of $u_{\pm}$ is a smooth bijection of $J$ onto $\R$,
and one thinks of $u_\pm$ as \lq\lq new coordinates" on $J$, adapted to
$f$.

 \medskip
It is possible to continue this process, to develop moduli
for more Takens cases, and for diffeomorphisms that are
flat at the ends. However
the flat (and semi-flat) cases offer enormous variety, and condition (E)
as it stands seems the simplest way to express the obstruction
to smooth conjugacy, given $C^1$ conjugacy. The modulus conditions
are computable in principle, but the computations are massive.

\section{Conjugacy in $\Diffeo^+_B$}\label{section-Diffeo+B}
Throughout this section, $I$ will be one of $\R$,
$[-1,+\infty)$, or $[-1,1]$. (Each
interval  with nonempty interior is
diffeomorphic to one of these, and it will be convenient
to have $0$ in the interior of $I$.)

$E$ will be a fixed closed
nonempty subset of $I$, containing any ends that $I$ has,
and $B$ will be the boundary of $E$.

\subsection{Proof of Lemma \ref{lemma-Diffeo+B-flat}}
The end $p$ in question is an accumulation point of $E$,
and hence $f-x$ and $g-x$ are flat there.
By assumption, there exists some $\phi_0\in\Diffeo^+_B$
that
conjugates
$g$ to $f$ on $I$. Since $\phi_0$ fixes each point
of $B$, $\phi_0-x$ is flat at $p$.
If $\phi$ is another map that conjugates $f$ to $g$ on $I$,
then $\phi_0^{\circ-1}\circ\phi$ belongs to the centraliser
$C_f$ of $f$ in $\Diffeo^+_B$ .  Since
$f-x$ is flat at $p$, so is every element of $C_f$
by Proposition \ref{proposition-flat-centraliser}.  Thus
$\phi$ is the composition of two functions that fix $p$ and have
Taylor series $X$ there, and the result follows. \qed

\subsection{Proof of Theorem \ref{theorem-Diffeo+B}}
We need only prove the \lq\lq if" part.

Assume Conditions (M) and (F).

We define a variation of $\Conj(f,g,J)$, corrsponding
to conjugations that fix each point of $B$:
\Definition If $L$
is a connected component of $I\sim B'$,
we denote the
set of all maps that conjugate
$f$ to $g$ in $\Diffeo_B^+(\clos(L))$
by $\Conj_B(f,g;L)$.

\medskip
The assumption (M) tells us that each $\Conj_B(f,g;L)$
is nonempty.
We have to show that we can patch
together elements
of the various $\Conj_B(f,g;L)$
to get an element of
$\Diffeo^+(I)$.

We claim that for each $L$, we may choose
$h\in\Conj_B(f,g;L)$ with $h-x$ flat
at each end of $L$ that is not an end of $I$.

Let $p$ be an end of $L$ that is not an end of $I$. Then $p\in B'$.
There are three cases:

1. $p$ is an accumulation point of  and $L\cap B$.
Then since all elements of $\Conj_B(f,g;L)$
fix all points of $L\cap B$, they all have
$h-x$ flat at $p$.

2. $p$ is isolated in $L\cap B$, and $p$
is an end of some component $J\subset L$
of $I\sim E$
(note that $p$ must be a limit
of points of $B$ that lie on the
other side of $p$ from $J$).
Then Condition (F) tells us that
all elements of $\Conj(f,g;J)$ have $h-x$
flat at $p$. Thus all elements of $\Conj(f,g;L)$
also have $h-x$ flat at $p$.

3.  $p$ is isolated in $L\cap B$,
and $p$ is an end of some component
$J\subset L$ of the interior of $E$.
Then given any $h\in \Conj(f,g;L)$,
we may modify it on $J$ in any way at all
(provided it remains a diffeomorphism
of $J$ onto itself) without disturbing
the conjugacy, because $f(x)=g(x)=x$ on $J$.
Thus we can modify it
to make $h-x$ flat at $p$.

\medskip
So the claim holds. So if we choose
$h$ on each $L$ to have $h-x$ flat at
each end in the interior of $I$, then they automatically fit
together to make
the desired conjugation.
\qed

\section{Reducing from $\Diffeo(I)$ to $\Diffeo^+(I)$}\label{section-Diffeo}

In this section we discuss the reduction of the conjugacy
problem in  the full diffeomorphism group
to the conjugacy problem in the subgroup
of direction-preserving maps.

There is no issue for half-open intervals, since the two groups
coincide, so it suffices to consider the two
cases $I=\R$ and $I=[-1,1]$, which represent all other
intervals up to diffeomorphism. (It is convenient to use
representatives that are invariant under $-$.)

\subsection{Reducing to conjugation by elements of $\Diffeo^+$}
The first (simple)
proposition allows us to restrict attention to
conjugation using $h\in\Diffeo^+(I)$.
\begin{proposition}\label{proposition-reduce-1}\label{general}
Let $I=\R$ or $[-1,1]$. Let $f,g\in\Diffeo(I)$. Then the following
two conditions are equivalent:\\
(1) There exists $h\in\Diffeo(I)$ such that $f=g^h$.
\\
(2) There exists $h\in\Diffeo^+(I)$ such that $f=g^h$ or
$-\circ f\circ-=g^h$.
\end{proposition}
\Proof
If (1) holds, and $\deg h=-1$, then
$-\circ f\circ - = g^k$, with
$$ k(x) = h(-x).$$
The rest is obvious.
\qed

\subsection{Reducing to conjugation of elements of $\Diffeo^+$}
The degree of a diffeomorphism is a conjugacy
invariant, so to complete the reduction of the
conjugacy problem in $\Diffeo$ to the problem
in $\Diffeo^+$, it suffices to deal with the
the case when $\deg f=\deg g=-1$ and
$\deg h=+1$.

\smallskip
Note that $\fix(f)$ and $\fix(g)$ are singletons, and lie in int$(I)$.
If $f=g^h$,
then $h(\fix(f))=\fix(g)$,
and (since $\Diffeo^+$ acts transitively on int$(I)$)
we may thus, without loss in generality, suppose
that $f(0)=g(0)=h(0)=0$.

If $f=g^h$, then we also have
$f^{\circ 2} =
(g^{\circ 2})^h$,
$f^{-1}=(g^{-1})^h$,
and $f^{\circ 2}\in\Diffeo^+$.

We have the following reduction:
\begin{theorem}\label{theorem-reduction}
Let $I=\R$ or $[-1,1]$.
Suppose $f,g \in
\Diffeo^-$, fixing $0$. Then the following two condition are
equivalent:
\begin{enumerate}
    \item $f=g^h$ for some $h\in \Diffeo^+$.
    \item \begin{enumerate}
        \item There exists $h_1\in \Diffeo^+_0$ such that $f^{\circ 2}
=(g^{\circ 2})^{h_1}$;

        and
        \item Letting $g_1=g^{h_1}$, there exists $h_2\in
        \Diffeo^+$, commuting with $f^{\circ 2}$ and fixing $0$,
        such that $T_0 f=(T_0 g_1)^{T_0 h_2}$.
    \end{enumerate}
\end{enumerate}
\end{theorem}
See \cite{OR}, Proposition 2.1 for details.

\subsection{Making the conditions explicit}
To complete the project of reducing
conjugation in $\Diffeo$ to conjugation in $\Diffeo^+$,
we have to find an effective way to check condition $2(b)$.
In other words, we have to replace the nonconstructive
\lq\lq there exists $h_2\in\Diffeo^+$"
by some condition that can be checked algorithmically.
This is achieved by the following:

\begin{theorem}\label{theorem-last}
Let $I=\R$ or $[-1,1]$.
Suppose that $f,g\in\Diffeo^-$ both fix $0$,
and have $f^{\circ 2}
=g^{\circ 2}$.
Then
there exists $h\in\Diffeo^+$, commuting with $f^{\circ 2}$,
such that $T_0f= (T_0g)^{T_0h}$
if and only if one of the following holds:
\begin{enumerate}
\item $(T_0f)^{\circ 2}\not=X$;
\item $0$ is an interior point of $\fix(f^{\circ 2})$;
\item $(T_0f)^{\circ 2} = X$, $0$ is a boundary point
of $\fix(f^{\circ 2})$, and $T_0f=T_0g$.
\end{enumerate}
\end{theorem}

Note that the conditions 1-3 are mutually-exclusive.
We record a couple of corollaries:

\begin{corollary}
Suppose $f,g \in
\Diffeo^-$, fixing $0$, and suppose $(T_0f)^{\circ2}\not=X$
or $0\in$int$\fix(f)$. Then
    $f=g^h$ for some $h\in \Diffeo^+$
if and only if
    $f^{\circ2}=(g^{\circ2})^h$ for some $h\in \Diffeo^+$.
\end{corollary}

In case $(T_0f)^{\circ2}\not=X$, any $h$ that conjugates $f^{\circ2}$
to $g^{\circ2}$ will also conjugate $f$ to $g$. In the other case
covered by this corollary, it is usually necessary to modify
$h$ near $0$.

\begin{corollary}
Suppose $f,g \in
\Diffeo^-$, fixing $0$, and suppose
$(T_0f)^{\circ2}=X$
and $0\in$bdy$\fix(f)$. Then
    $f=g^h$ for some $h\in \Diffeo^+$
if and only if
    $f^{\circ2}=(g^{\circ2})^h$ for some $h\in \Diffeo^+$ and
$T_0f=T_0g$.
\end{corollary}

The last corollary covers the case where $0$ is
isolated in $\fix(f^{\circ2})$ and $T_0f$ is involutive,
as well as the case where $0$ is both an accumulation
point and a
boundary point of $\fix(f)$

The detailed proofs may be found in \cite{OR}. They use
results about conjugacy and reversibility
for formal power series, together with Kopell's
results about centralisers.

\section{Further Examples and Remarks}\label{section-examples}
\subsection{Reduction to $\fix(f)=\fix(g)$}
\label{examples-subsection-Diffeo^+}
In relation to the reduction of
Subsection \ref{subsection-Diffeo+},
it is not true that each map $h$ conjugating $g$ to $f$ may be factored
as {\em any} smooth map that maps $\fix(f)$ onto $\fix(g)$, followed by a
smooth map fixing $\fix(f)$.

\begin{example}
\end{example}
Take, for instance $f(x)= x+ \sin(x)/10$, and $g(x)=x-\sin(x)/10$.
Both fix precisely $\pi\Z$. They are conjugated by $h:x\mapsto
x-\pi$.  The identity map $h_1$ maps $\fix(f)$ onto $\fix(g)$, but
no map fixing $\fix(f)$ conjugates $f$ to $g$, since the multipliers
are wrong.

In general, in searching
for a factor  $h_1$ as in Proposition \ref{proposition-reduction-2},
we may start by classifying the points $p$ of $\bdy(\fix(f))$ (and $\bdy(\fix(g))$)
according to the conjugacy class of $T_pf$ (or $T_p(g)$).
This produces two classifications, the $f$-classification of $\fix(f)$, and
the $g$-classification of $\fix(g)$.  Only maps $h_1$ that
respect these classifications are eligible as potential factors.
Precisely speaking, the eligible maps $h_1$ must be such that
$T_pf$ and $T_{h_1(p)}g$ are conjugate Taylor series,
for each $p\in\bdy(\fix(f))$.

\begin{example}
\end{example}
For instance, if we modified the above example by taking
$$ f(x) = x + {\sin x\over 1+x^4},$$
then there is no eligible map at all, so $f$ and $g$ are not conjugate.

\begin{example}
\end{example}
If we modified $g$ as well, taking
$$ g(x) = x + {\sin 2x\over 2+ 8x^2} \quad (= \half f(2x) ),$$
then the only eligible $h_1$ are those that have
$h_1(x)=2x$ on $\fix(f)=\pi\Z$.

This prompts the question, whether, assuming
the maps $f$ and $g$ are conjugate,
 {\em every} $h_1\in\Diffeo^+$
that respects this Taylor-series classification at the boundary
points will serve as a factor of the kind
referred to in Proposition \ref{proposition-reduction-2}.
That would be very convenient, as it would characterise
the diffeomorphisms $h_1$ that we need to find.
Unfortunately, the answer is no:

\begin{example}
\end{example}
Take any $f\in\Diffeo^+(\R)$ that fixes precisely $\Z$, has
$f-x$ flat at each integer, and is such that the functions
on $[0,1]$ defined by
$x\mapsto f(x+n)-n$ ($n\in\Z$) represent distinct
conjugacy classes of $\Diffeo^+([0,1])$.
Take $g(x)=1+f(x-1)$.
Then $f$ and $g$ are conjugate, but the map
$h_1(x)=x+2$ won't do as a factor of the required
kind, because no map that fixes $\Z$ will
conjugate $g$ to $x\mapsto 2+f(x-2)$.

\medskip
It is not essential to use a function $f$ that is flat
on the boundary to give an example of this kind.
We know that in the non-flat case, Condition (T)
is not enough to characterise conjugacy in
$\Diffeo^+(I)$, for compact $I$, so we can modify
the example to produce the same end result
without having $f-x$ flat at all.
The point is that once we have a $C^\infty$ diffeomorphism on each
interval $[n,n+1]$ and the two available Taylor series agree
at each $n$, then they patch together to make a
global diffeomorphism.

\medskip
So there is a substantial problem, from the constructive point-of-view,
concerning how to search for suitable $h_1$.

However, we know, from Subsection \ref{subsection-modulus}
that part (3) of Proposition \ref{proposition-reduction-2}
can only work if $f_1$ and $g$ have the same $J$-modulus,
for each component $J$.  This provides a fine filter,
to cut down the search, because given the Taylor
series at one end, there is at most a one-parameter
coset of diffeomorphisms of $J$ that conjugate $f_1$ to $g$
on $J$.
Generically, the coset is discrete.

\subsection{Finite $\fix(f)$}\label{examples-subsection-finite}
We conclude with a summary of our conclusions about
 the conjugacy problem in the
special case when $\fix(f)=\fix(g)=E=B$ is a finite set
of points $p_1<p_2<\cdots p_k$,
so that there is only
only one way to map $E=B$ to itself, preserving order.

The case when all the points are hyperbolic is classical, and
has been discussed previously by Belitsky \cite{B}.
We include this case in the discussion, for completeness.
The first necessary condition is that Condition (T)
holds at each of the fixed points, i.e. that the Taylor series of
$f$ and $g$ be conjugate. At hyperbolic points for $f$ or $g$, this amounts
to the identity of the multipliers, and at the remaining points
$p$ at which $f-x$ is not flat (\lq\lq Takens points"),
it is determined by examining
the coefficients of $f$ and $g$ up as far as the term in
$X^{2p+1}$, where $T_pf-X$ vanishes to order $p$, but not
to order $p+1$.  At the points where $f-x$ is flat, the condition is
automatic.

Next, we need the sign condition, that the graphs of $f$ and $g$ lie
on the same side of the diagonal on each interval complementary
to the fixed-point set.  This condition follows automatically
from (T) at the hyperbolic and Takens points.

Next, we need condition (P), the convergence of the products (1).
This is automatic at the hyperbolic and Takens
points, but imposes restrictions to the right and left of the points where
$f-x$ is flat.

Next, we need condition (E), to the effect that the $C^1$
conjugacies that now exist include some that are $C^\infty$
when restricted to each of the half-open intervals
$(-\infty,p_1]$, $[p_1,p_2)$, $(p_1,p_2]$, $[p_2,p_3)$,
$(p_2,p_3]$,$\ldots$,$[p_k,+\infty)$\footnote{This much is automatic at the
hyperbolic and Takens points.}, and that on each
of the compact intervals $[p_1,p_2]$,$\ldots$,
$[p_{k-1},p_k]$ there is at least one of these
conjugacies that is smooth to both ends.
Note that this means that the maps $f$ and $g$
share the same Robbin invariant.

At this stage, we have nonempty
cosets $\Conj(J)=\Conj(f,g;J)$ of maps that conjugate
smoothly on each clos$J$. Next we need
condition (M),
that we can match some  Taylor series
from $\Conj(J)$ and $\Conj(J')$ whenever $J$ and $J'$ are
adjacent components.  This may still not be enough
to make $f$ and $g$ conjugate.

We distinguish the hyperbolic and Takens points
from the points where $f-x$ is flat.  At the latter,
all smooth $J$-conjugations share the same
Taylor series, as do all $J'$-conjugations, so if
any series from $\Conj(J)$ concides with a series
from $\Conj(J')$, then all do, so we can stop
worrying about these fixed points.

Let $q_1<\cdots<q_r$ be the remaining fixed points, the
ones at which $f-x$ is not flat.  Write $L_i=(-\infty,q_i]$,
and $J_i=[q_i,q_{i+1}]$, ($i<r-1$) and $J_r=[q_r,+\infty)$.

If $r=1$, we are done; $f$ and $g$ are smoothly conjugate.
Otherwise $\Conj(L_2)$ is already nonempty. We have to assume
that $\Conj(L_2)\cap \Conj(J_2) = \Conj(L_3)$ is nonempty; otherwise
$f$ and $g$ are not conjugate. Checking this condition is
a matter of comparing the set of multipliers (at hyperbolic points)
or the set of $2p+1$-st order Taylor polynomials
(at Takens points). Each of these sets is a coset of a group.
In the hyperbolic case, unless one of the sets
is the full multiplicative group $(0,+\infty)$,
we are comparing two sets of the form
$$\{ \alpha\lambda^n: n\in \Z\}\hbox{ and }
\{\beta\mu^n: n\in\Z\},$$
(i.e. two cosets of discrete subgroups of
the $(0,+\infty)$).
In the Takens case, once we conjugate the series to canonical
form, we are comparing the coefficients of $X^{p+1}$, which
are two cosets of the additive group $\R$.

Continuing, we get
a decreasing sequence $\Conj(L_4)$,$\ldots$,$\Conj(L_r)$,
and if at any stage it is empty, there is no conjugacy.

If $\Conj(L_r)$ is nonempty, the last step is to see
whether $\Conj(L_r)\cap \Conj(J_r)$ is nonempty. If this last condition
holds, then there is a smooth conjugacy between f and g, and
otherwise not.

\subsection*{Acknowledgements}
The authors are most grateful to \'Etienne Ghys and Ian Short
for useful advice, and to John Mather for making
available his unpublished work.

\bigskip
\noindent
e-mail:\\
anthonyg.ofarrell@gmail.com\\
maria@chalmers.se

\end{document}